\documentclass[%
 aip,
 amsmath,amssymb,
reprint,%
]{revtex4-1}

\usepackage{graphicx,subfigure,varwidth}
\usepackage{dcolumn}
\usepackage{bm}
\usepackage[utf8]{inputenc}
\usepackage[T1]{fontenc}
\usepackage{mathptmx}
\usepackage{etoolbox}
\usepackage{array}
\usepackage[
            citecolor=blue]{hyperref}

\def\rmd{\mathrm{d}}
\def\rme{\mathrm{e}}

\def\mcR{{R}}

\makeatletter
\def\@email#1#2{%
 \endgroup
 \patchcmd{\titleblock@produce}
  {\frontmatter@RRAPformat}
  {\frontmatter@RRAPformat{\produce@RRAP{*#1\href{mailto:#2}{#2}}}\frontmatter@RRAPformat}
  {}{}
}%
\makeatother
\begin{document}

\preprint{AIP/123-QED}
\title{Two-stage Fourth-order Gas Kinetic Solver-based Compact Subcell Finite Volume Method for Compressible Flows over Triangular Meshes}
\author{Chao Zhang}
\affiliation{ 
Institute of Applied Physics and Computational Mathematics, Beijing, 100088, China
}%
\author{Qibing Li}
\affiliation{ 
Department of Engineering Mechanics, Tsinghua University, Beijing 100084, China
}
\author{Peng Song}
\affiliation{ 
Institute of Applied Physics and Computational Mathematics, Beijing, 100088, China
}
\affiliation{HEDPS, Center for Applied Physics and Technology,  College of Engineering, Peking University, Beijing 100871, China.}
\author{Jiequan Li}
\affiliation{ 
Institute of Applied Physics and Computational Mathematics, Beijing, 100088, China
}
\affiliation{HEDPS, Center for Applied Physics and Technology,  College of Engineering, Peking University, Beijing 100871, China.}
\email{li\_jiequan@iapcm.ac.cn.}
 
\begin{abstract}
To meet the demand for complex geometries and high resolutions of small-scale 
flow structures, a two-stage fourth-order subcell finite volume (SCFV) method  combining the gas-kinetic solver (GKS) with subcell techniques for compressible flows over (unstructured) triangular meshes was developed to improve the compactness and efficiency. Compared to the fourth-order GKS-based traditional finite volume (FV) method, the proposed method realizes compactness effectively by subdividing each cell into a set of subcells or control volumes (CVs) and selecting only face-neighboring cells for high-order compact reconstruction. Because a set of CVs share a solution polynomial, the reconstruction is more efficient than that for traditional FV-GKS, where each CV needs to be separately reconstructed. Unlike in the single-stage third-order SCFV-GKS, both accuracy and efficiency are improved significantly by two-stage fourth-order temporal discretization, for which only a second-order gas distribution function is needed to simplify the construction of the flux function and reduce computational costs. For viscous flows, it is not necessary to compute the viscous term  with GKS. Compared to the fourth-stage Runge--Kutta method, one half of the stage is saved for achieving fourth-order time accuracy, which also helps to improve the efficiency. Therefore, a new high-order method with compactness, efficiency, and robustness is proposed by combining the SCFV method with the two-stage gas-kinetic flux. Several benchmark cases were tested to demonstrate the performance of the method in compressible flow simulations. 
\end{abstract}

\maketitle

\section{\label{sec:level1}Introduction}
Increasing engineering demands for accuracy and efficiency have catalyzed the development of high-order methods in the field of computational fluid dynamics. Complex geometries are often involved in engineering problems, necessitating the use of unstructured meshes, which are more flexible and more easily generated than structured meshes. 
Several flow problems require low dissipation and dispersion errors, such as large eddy simulation and direct numerical simulation of turbulence and problems in aeroacoustics.
Low-order methods are usually too dissipative to resolve small-scale structures. In contrast, high-order methods can potentially achieve increased accuracy with low computational cost. However, existing high-order methods are generally less robust and more complicated than low-order methods. Therefore, it is necessary to develop new high-order methods for unstructured meshes to improve robustness, simplicity, and efficiency.

This study was performed in the finite volume (FV) framework. Compared to finite difference (FD) methods, finite volume methods, and finite element (FE) methods are more suitable for unstructured meshes. FV methods are usually simpler and more robust than FE methods, especially for high-speed flows containing discontinuities (e.g., shocks). Many high-order FV methods were developed on unstructured meshes, such as the k-exact method \cite{TJBarth1990}, the essentially non-oscillatory (ENO) method \cite{RAbgrall1992}, and the weighted ENO (WENO) method \cite{OFriedrich1998}. Because only cell averages are available in common FV methods, solution polynomials with higher degrees of freedom (DoFs) are required to represent the flow structures, which may reduce the compactness, efficiency, and robustness. For example, the non-compactness resulting from a large stencil for reconstruction leads to problems with numerical dissipation, high-order numerical approximations of boundary conditions, reduced parallel efficiency, and caches missing. 

To achieve compactness for high-order FV methods, the compact least-square reconstruction and variational reconstruction \cite{QWang2016II, QWang2017} were developed, where a large linear equation system should be solved and implicit time-stepping methods are necessary. 
By subdividing each cell into a set of subcells, the spectral volume (SV) achieves compactness as well\cite{ZJWang2004}, and subcell-averaged solutions are updated under the FV framework. Compared to traditional FV methods, the SV method can achieve higher resolution with the same mesh size owing to the use of subcells, which helps to capture flow structures with small scales more accurately. 
However, the SV method requires the number of subcells to be equal to that of the unknowns in the solution polynomial. The complex subdivision of the main cells thus becomes a bottleneck, and no stable subdivision has yet been found for quadrilateral or hexahedral meshes. As an extension of the SV method, the subcell finite volume (SCFV) method \cite{JHPan201708,JHPAN2017} overcomes the difficulty of subdivision by combining the idea of the $\mathrm{P_NP_M}$ method, avoiding the restriction of the number of subcells in SV. The subcells from face-neighboring cells are involved in the stencil for reconstruction. As a result, the subdivision is much easier and more flexible. Theoretically, the SCFV method can achieve arbitrarily high-order accuracy with a compact stencil. 
In summary, the SCFV method has the following three advantages over existing high-order methods:

(i) Compact reconstruction is more easily implemented compared to traditional FV methods. Generally,
the arbitrary distribution and number of unstructured meshes make it difficult to choose and search neighboring cells to form the stencil. The SCFV method combines the features of both unstructured and structured meshes. First, a set of coarser unstructured meshes are generated, and are called main cells. Each main cell is further subdivided into a set of subcells that are stored in a structured way. In fact, these subcells are the real control volumes (CVs), and the subcell averages are updated in the FV framework. In this way, the complexity of choosing and searching neighbors is reduced greatly. As only face-neighboring main cells are considered, there are usually enough DOFs, i.e., subcells, for a compact reconstruction.

(ii) The reconstruction is efficient compared to traditional FV methods.
Assume that each main cell is subdivided into $n$ subcells. Similar to finite element (FE) methods, the $n$ subcells (CVs) are treated as inner DOFs, and a common solution polynomial is reconstructed for them. With the same number of CVs, reconstruction is conducted only $1/n$ times compared to traditional FV methods for which the reconstruction needs to be conducted for each CV separately. Moreover, the memory space for the corresponding coefficient matrix is also $1/n$ of that in traditional FV. Therefore, the reconstruction in the SCFV method is more efficient in terms of both computational cost and memory storage.

(iii) The subcell resolution can be well preserved compared to FE methods.
There are many FE or hybrid methods that also achieve compactness easily, such as the DG method \cite{BCockburn1998}, the CPR method \cite{Huynh2007,ZJWang2009}, and the $\mathrm{P_NP_M}$ method \cite{MDumbser2008}. However, the solution distribution needs to remain continuous inside each main cell, which results in the difficulty with respect to capturing shocks. It is difficult to fully maintain the advantage of inner DOFs, i.e., the subcell resolution. In contrast, for the SCFV method, discontinuities can even exist inside each main cell. The resolution and robustness for shock capturing can thus be enhanced effectively. Besides, the SCFV method can take larger CFL numbers, and the formulation is much simpler compared to FE methods.

The SCFV method was originally developed on quadrilateral meshes to solve the Euler equations \cite{JHPan201708,JHPAN2017}. Then, it was further developed on triangular meshes and extended to the NS equations by combining the single-stage third-order gas-kinetic flux solver (GKS) \cite{CZhang2021}. Although other solvers, such as the GRP solver, \cite{MBenArtzi2006,MBenArtzi2007} can also be adopted, this study still uses the GKS solver but with second-order accuracy \cite{KXu2001,QBLi2006} for the flux construction because the inviscid and viscous flux are coupled and computed simultaneously \cite{LZhang201931,NZhan2021,LMYang202133}. Besides, a series of unified gas-kinetic scheme  (UGKS) has been developed  for  the  entire  flow  regime  from  the  continuum  flow  to  the 
highly rarefied flow\cite{KXu2010229,YZhang201931,YChen202032,MZhong202133}.
During the past decade, high-order GKS has been developed systematically. By taking a second-order Taylor expansion, a third-order multi-dimensional GKS was developed within a single stage \cite{QBLi2010}. The time-dependent gas distribution function can even be used to evolve the solutions at cell interfaces, which provides additional DoFs to help construct a compact third-order FV-GKS \cite{LPan2016318,XJi2020410}. Other extensions of single-stage third-order GKS can also be found for the DG method \cite{XDRen2015} and the CPR method \cite{CZhang2018}. As mentioned above, a single-stage third-order SCFV-GKS has been proposed recently \cite{CZhang2021}. However, when fourth-order accuracy is considered, the single-stage flux function is quite complicated, and the computational cost increases significantly \cite{NLiu2014}; the robustness may also be weakened. Fortunately, the two-stage fourth-order time-stepping method \cite{JLi201638} provides an efficient way for Lax-Wendroff-type solvers to achieve fourth-order time accuracy. Only the second-order gas-kinetic flux solver is needed, which is much simpler and more robust \cite{FXZhao2019,LPan2021,FXZhao202133}. The computational cost of the flux evaluation can be reduced effectively compared to the single-stage high-order GKS\cite{LPan2016II}.

Hence, this study developed a two-stage temporal-spatial fourth-order SCFV-GKS to solve compressible flows over triangular meshes, with the aim of enhancing the compactness, efficiency, and robustness. Compared to the single-stage SCFV-GKS, both accuracy and efficiency are improved significantly. By adopting the two-stage temporal discretization, the flux evaluation under the SCFV framework can be simplified greatly, making it easier to implement. The robustness is also enhanced because of the use of the second-order gas-kinetic flux. These improvements make the SCFV method more promising for engineering applications.

This paper is organized as follows. Section 2 presents the compact fourth-order reconstruction based on the SCFV method and the two-stage gas-kinetic flux evolution. Section 3 presents numerical tests with several benchmark cases to demonstrate the high accuracy, efficiency, and robustness of the current method for simulating compressible flows. The last section presents the conclusions.

\section{Gas kinetic solver-based subcell finite volume method}

In this section, we propose a GKS-based SCFV method for compressible fluid flows. The primary steps are the compact subcell technique and the two-stage fourth-order temporal advancing based on the second-order GKS. 

\subsection{Subcell finite volume method}

\begin{figure}[!htb]
  \centering
  \includegraphics[scale=0.13]{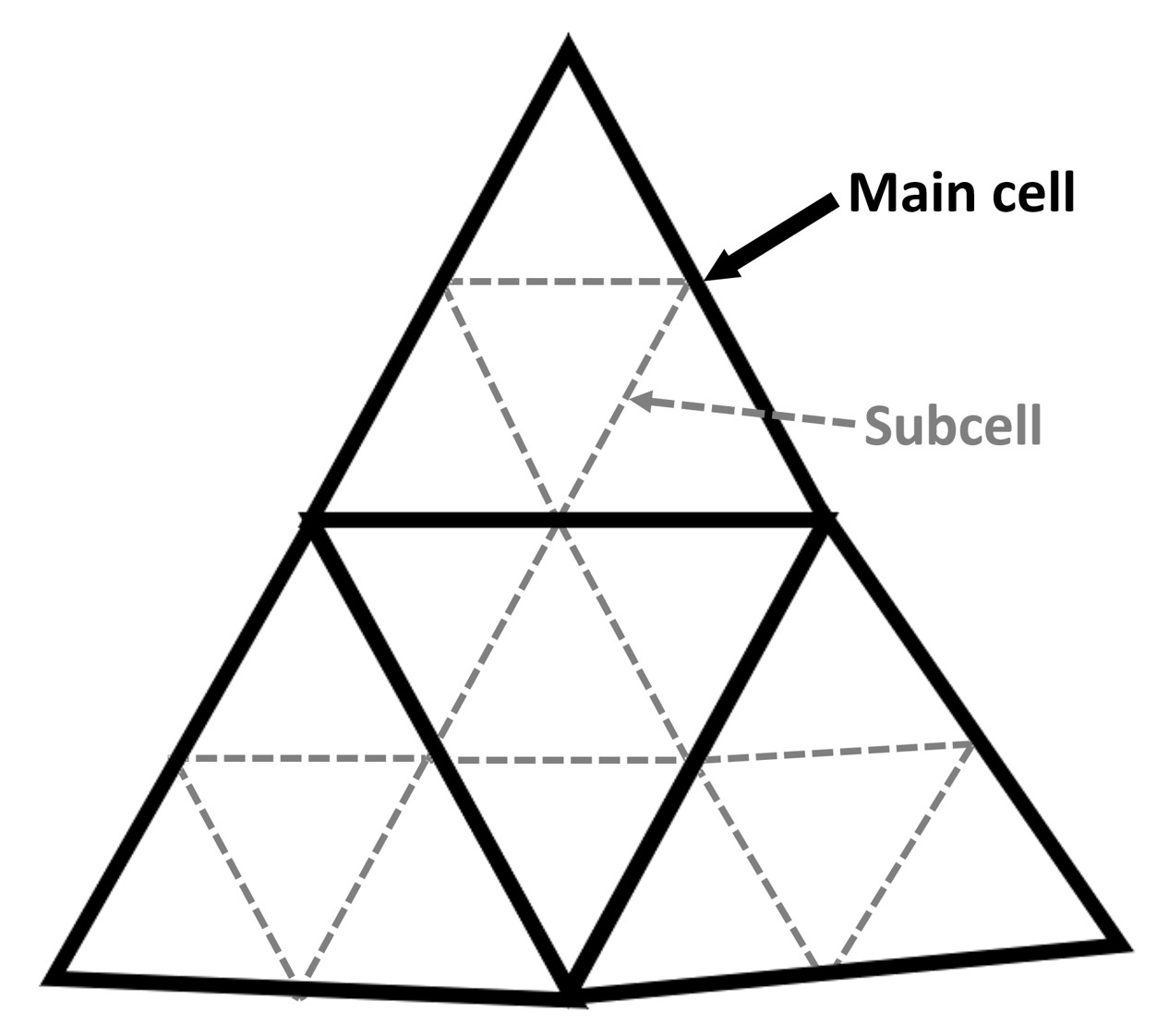}\\
  \caption{Subdivision of main cells.}
  \label{fig:SCFV_Stencil}
\end{figure}

First, we briefly review the SCFV method. Consider the two-dimensional (2D) conservation law
\begin{equation}\label{eq_ConLaw}
    \frac{\partial\bm{Q}}{\partial{t}}+\nabla\cdot\bm{F} =0,
\end{equation}
where $\bm{Q}=(\rho,\rho\bm{U},\rho E)^T$ are the conservative variables, where $\bm{U}=(U,V)$ are the macroscopic velocities and $\bm{F}$ is the flux vector. The computational domain is divided into $N$ non-overlapping triangular cells $\{\Omega_i\}$, which are also referred to as main cells. As shown in Figure \ref{fig:SCFV_Stencil}, each main cell is uniformly subdivided into four similar subcells $\{\Omega_{i,j}\}$. 
In the SCFV method, subcell-averaged solutions are stored and updated under the FV framework. More clearly, a subcell can also be called a CV when compared to that for traditional FV methods. The primary task is to reconstruct solution polynomials over each main cell. Then, the numerical flux at subcell interfaces can be determined and used to update subcell averages. The reconstruction is implemented for each component of conserved variables $Q$. To achieve fourth-order space accuracy, $Q$ over each main cell ${\Omega}_i$ is approximated by a solution polynomial $P_i$ of degree $k = 3$, where $(k+2)(k+1)/2=10$ unknown coefficients need to be determined. As shown in Figure \ref{fig:SCFV_Stencil}, by involving face-neighboring cells, there are 16 subcells available to sufficiently determine the unknowns. 

\begin{figure*}[!htb]
  \centering
  \includegraphics[scale=0.23]{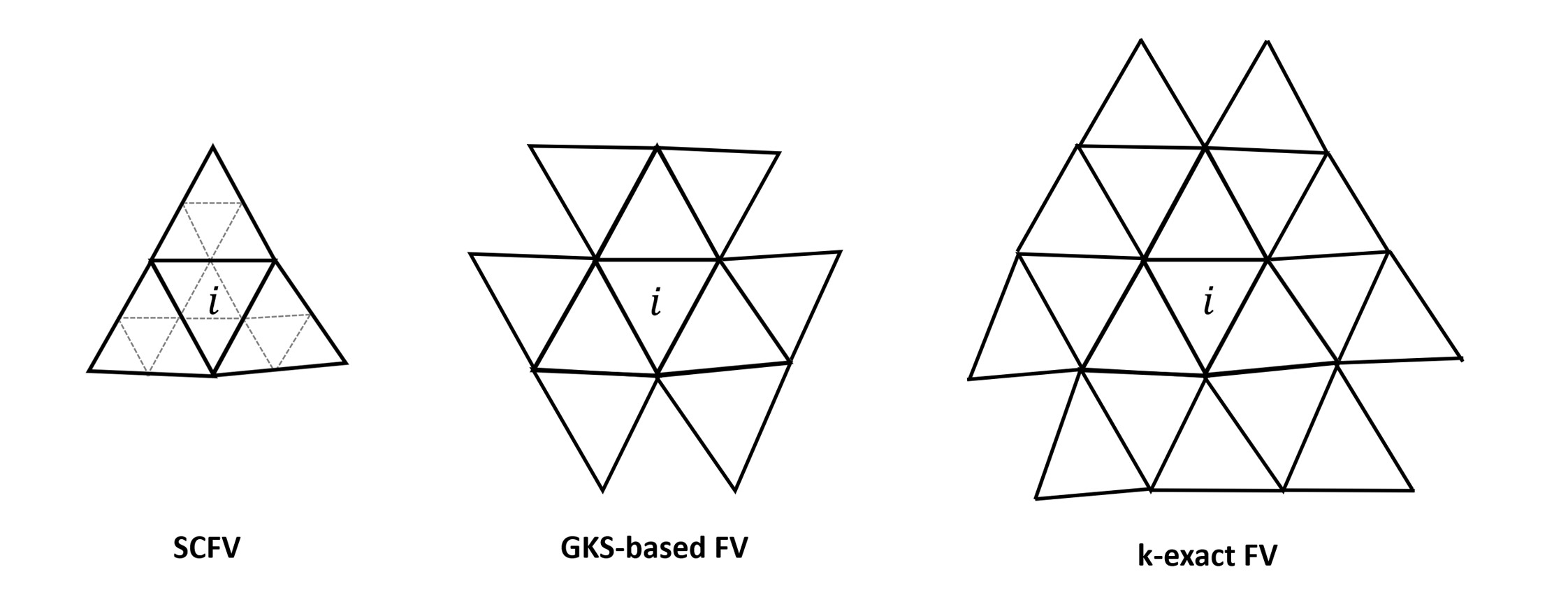}\\
  \caption{Comparison of the stencils for three types of fourth-order reconstructions: SCV, GKS-based FV, and k-exact FV.}
  \label{fig:SCFV_Stencil_compare}
\end{figure*}

In comparison, Figure \ref{fig:SCFV_Stencil_compare} also shows the stencil used in a GKS-based FV method and the classical k-exact FV \cite{FXZhao2021Preprint}. Note that for the GKS-based FV, cell-averaged slopes are updated and participate in the reconstruction, which helps to avoid the large stencil for the k-exact FV. Nonetheless, with only face neighbors, it is still not enough for a fourth-order reconstruction, which results in the enlargement of the stencil. In contrast, the stencil of SCFV is most compact. 

\begin{figure*}
  \centering
  \includegraphics[scale=0.3]{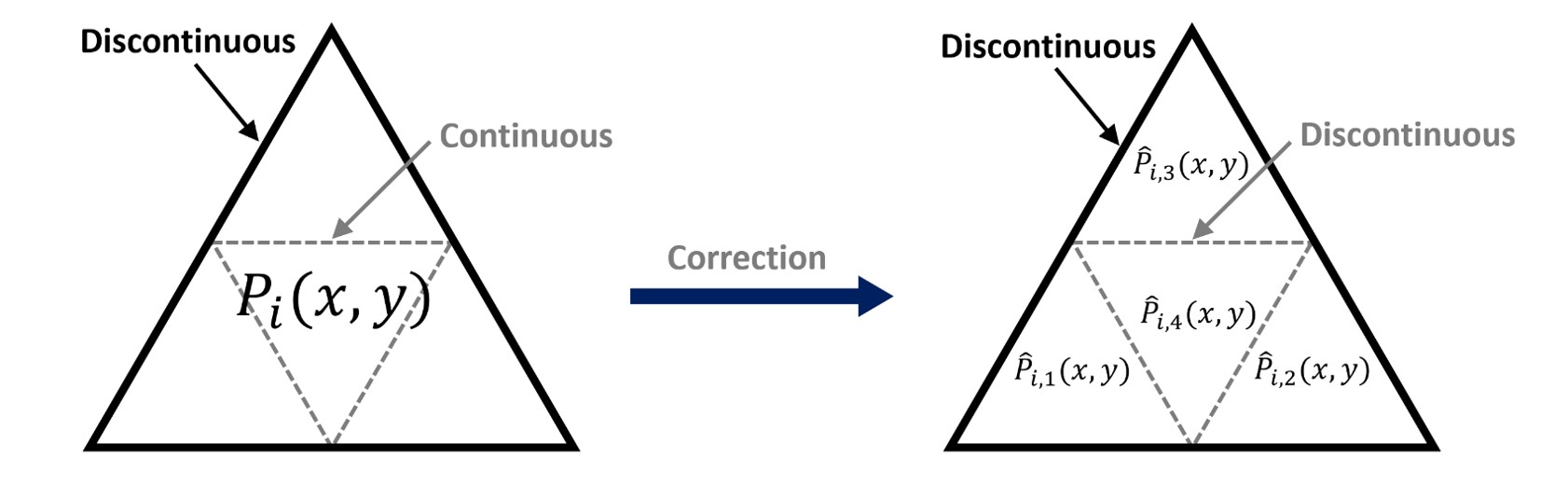}\\
  \caption{Correction for conserving subcell averages.}
  \label{fig:SCFV_Correction}
\end{figure*}

 For the convenience of further description, all subcells involved in the stencil are denoted as ${\Omega}_{i,st}~(st=1,\cdots,16)$. Besides, ${\Omega}_{i,j}~(j=1,\cdots,4)$ is specifically used to indicate subcells from the target main cell ${\Omega}_{i}$. The solution polynomial on ${\Omega}_{i}$ is expressed as
\begin{equation}\label{eq_Q_Poly}
\begin{split}
   P_i(x,y)=&\bar{Q}_i
   				+Q_{i,x}\left(\Delta{x}-\widehat{x}\right)
   				+Q_{i,y}\left(\Delta{y}-\widehat{y}\right)\\
   			  &+\frac12Q_{i,xx}\left(\Delta{x}^2-\widehat{x^2}\right)
   				+Q_{i,xy}\left(\Delta{x}\Delta{y}-\widehat{xy}\right)\\
   				&+\frac12Q_{i,yy}\left(\Delta{y}^2-\widehat{y^2}\right)
   				+\frac16Q_{i,xxx}\left(\Delta{x}^3-\widehat{x^3}\right)\\
   				&+\frac12Q_{i,xxy}\left(\Delta{x}^2\Delta{y}-\widehat{x^2y}\right)+\frac12Q_{i,xyy}\left(\Delta{x}\Delta{y}^2-\widehat{xy^2}\right)\\
   				&
   				+\frac16Q_{i,yyy}\left(\Delta{y}^3-\widehat{y^3}\right),
\end{split}
\end{equation}
where $\Delta{x}=x-x_{i,c},~\Delta{y}=y-y_{i,c}$,  $(x_{i,c},y_{i,c})$ indicates the centroid of ${\Omega}_{i}$. 
The zero-mean basis functions are adopted so that the averaged solution $\bar{Q}_i$ can be conserved automatically on ${\Omega}_{i}$. Thus, there are $9$ unknowns $Q_{i,x}, Q_{i,y},\cdots,Q_{i,yyy}$ to be determined for each $\Omega_i$. The definition of $\widehat{x^\alpha y^\beta}$ is
\begin{equation}\label{eq_xy_int}
  \widehat{x^\alpha y^\beta}=\frac{1}{|{\Omega}_{i}|} \int_{\Omega_i}
          \Delta{x}^\alpha \Delta{y}^\beta\mathrm{d}\Omega,
\end{equation}
where $\alpha,\beta=0, \cdots, 3$ and satisfy $0\leq\alpha+\beta\leq3$. $|{\Omega}_{i}|$ is the volume of ${\Omega}_{i}$. Integrating $P_i(x,y)$ over ${\Omega}_{i,st}$ gives
\begin{equation}\label{eq_eq_set}
    Q_{i,x}\overline{x}^{i,st}+Q_{i,y}\overline{y}^{i,st}+...+\frac16Q_{i,yyy}\overline{y^3}^{i,st}=\bar{Q}_{i,st}-\bar{Q}_i,
\end{equation}
where $st=1,\cdots, 16$, 
\begin{equation}\label{eq_base_int}
    \overline{x^\alpha y^\beta}^{i,st}=\frac{1}{|\Omega_{i,st}|} \int_{\Omega_{i,st}}\left(\Delta{x}^\alpha\Delta{y}^\beta-\widehat{x^\alpha y^\beta}\right)\mathrm{d}\Omega,
\end{equation}
and $\bar{Q}_{i,st}$ is the averaged solution on the subcell $\Omega_{i,st}$. Note that because $\bar{Q}_i$ has been used in Eq.(\ref{eq_Q_Poly}) and it has the relation with $\bar{Q}_{i,j}$ as
\begin{equation}\label{eq_Qij_Qi}
    \sum_{j=1}^4 \bar{Q}_{i,j} |\Omega_{i,j}|=\bar{Q}_i |\Omega_i|,
\end{equation}
 the target main cell ${\Omega}_{i}$ can only provide three more DOFs for reconstruction besides $\bar{Q}_i$. Here, the central subcell is excluded from the stencil. Therefore, 15 equations are left in Eq.(\ref{eq_eq_set}), which forms an over-determined system. The unknowns are solved by the weighted least-square technique
\begin{equation}\label{eq_WLS}
  \mathrm{min}\left[\sum_{st} w_{i,st}\left(\frac{1}{|\Omega_{i,st}|} \int_{\Omega_{i,st}}P_i(x,y)\mathrm{d}\Omega-\bar{Q}_{i,st}\right)^2\right],
\end{equation} 
where the weight $w_{i,st}=1/d_{i,st}$, $d_{i,st}$ indicates the distance between the centroid of $\Omega_{i,st}$ and $\Omega_i$. The resulting solution polynomial $P_i(x,y)$ is shared by the four subcells ${\Omega}_{i,j}~(j=1,\cdots, 4)$. However, the averaged solutions of the four subcells are generally not conserved directly by the above weighted least-square reconstruction, i.e.,
\begin{equation}\label{eq_P_neq_Q}
  \frac{1}{|\Omega_{i,j}|} \int_{\Omega_{i,j}}P_i(x,y)\mathrm{d}\Omega\neq\bar{Q}_{i,j},~j=1,\cdots, 4.
\end{equation}
Studies showed that \cite{JHPan201708,JHPAN2017}, it is necessary to conserve the subcell averages on the four subcells, otherwise, the scheme may not be stable. A simple correction is adopted by shifting $P_i(x,y)$ on each subcell $\Omega_{i,j}$, i.e., $\hat{P}_{i,j}(x,y)=P_i(x,y)-\bar{P}_{i,j}+\bar{Q}_{i,j}$ where
\begin{equation}\label{eq_P_ave}
  \bar{P}_{i,j}=\frac{1}{|\Omega_{i,j}|} \int_{\Omega_{i,j}}P_i(x,y)\mathrm{d}\Omega.
\end{equation}
Then, $\hat{P}_{i,j}(x,y)$ is used for the flux evaluation.
The correction introduces jumps artificially across interfaces between subcells, as shown in Figure \ref{fig:SCFV_Correction}. The computational cost for the flux evaluation increases at interfaces between subcells. 
Some remarks are made below. 
Because the second-order gas distribution function is much simpler than the third-order one, and the computational cost for numerical flux construction has been reduced significantly, we simply take the weighted least-square (WLS) technique to determine the unknowns, rather than the constrained least-square (CLS) technique adopted in the single-stage SCFV-GKS \cite{CZhang2021}. Both the techniques have advantages and disadvantages,. With the CLS technique, the subcell averages can be conserved directly and the distribution inside each main cell can remain continuous so that the computational cost of flux is reduced. Thus, it introduces more benefits for the third-order gas-kinetic flux. However, the technique itself is more complicated and time-consuming, and it may reduce the robustness of SCFV. In contrast, the WLS technique is simpler and more efficient at the expense of increasing the computational cost of flux because of the correction. However, the second-order gas-kinetic flux is cheaper, and by introducing jumps, the robustness of SCFV is enhanced, and the resolution for discontinuities can be improved, especially for shock capturing.

For smooth flows, the solution polynomial obtained by the above reconstruction can be used directly for flux evaluation. However, for shock capturing, it may lead to numerical oscillations near shock waves. It is necessary to apply an effective limiting procedure or limiter to the solution polynomial. To reduce the computational cost, a shock detector \cite{WLi2012} is adopted to mark troubled cells near shock waves. For other unmarked cells, no limiter is needed and the high accuracy in smooth flow regions can be maintained. For troubled cells, a limiting procedure based on hierarchical reconstruction (or HR limiter) is implemented to obtain a new solution polynomial, which is able to suppress oscillations effectively\cite{ZLXu2009}. The limiter can keep the designed order of accuracy, and only face-neighboring cells are needed, which remains the compactness of the current scheme. For the details of the limiting procedure, we refer to Ref.\onlinecite{CZhang2021}.

\subsection{Two-stage gas-kinetic flux evolution}

By integrating Eq.(\ref{eq_ConLaw}) over each subcell $\Omega_{i,j}$, the semi-discrete form of the SCFV framework can be obtained as follows.
\begin{equation}\label{eq_semi_SCFV}
\begin{split}
    &\frac{\partial\bar{\bm{Q}}_{i,j}}{\partial{t}}=\mcR_{i,j}({\bm{F}}), \\
    &\mcR_{i,j}({\bm{F}})=-\frac{1}{|\Omega_{i,j}|}\oint_{\partial \Omega_{i,j}}\left({\bm{F}}\cdot\bm{n}\right)\mathrm{d}\Gamma,
\end{split}
\end{equation}
where $|\Omega_{i,j}|$ is the area of $\Omega_{i,j}$, $\partial \Omega_{i,j}$ is the interfaces surrounding $\Omega_{i,j}$, and $\bm{n}$ is the outward unit normal vector. To maintain fourth-order space accuracy, the line integral on the right-hand side is discretized using the Gaussian quadrature rule with two points
\begin{equation}\label{eq_F_gauss}
   \oint_{\partial \Omega_{i,j}}\left({\bm{F}}\cdot\bm{n}\right)\mathrm{d}\Gamma=
   \sum_{s\in\partial{\Omega}_{i,j}}\sum_{l=1}^{2}\omega_l({\bm{F}}\cdot\bm{n})_{s,l}|\Gamma|_s.
\end{equation}
where $|\Gamma|_s$ is the length of subcell interfaces, and $\omega_1=\omega_2=1/2$ are the
quadrature weights of the two Gaussian points. These Gaussian points are also referred to as flux points for clarity. For the temporal discretization of Eq.(\ref{eq_semi_SCFV}), by fully exploiting the time-evolving flux function provided by the second-order gas-kinetic flux solver, fourth-order time accuracy can be achieved more efficiently with the help of the two-stage temporal discretization \cite{JLi201638}, which can be expressed with the current notations as
\begin{equation}\label{eq_S2O4}
\begin{split}
&\bar{\bm{Q}}^*=\bar{\bm{Q}}^n+\frac12 \Delta t \mcR({\bm{F}}^n)+\frac18 \Delta t^2 \frac{\partial \mcR({\bm{F}}^n)}{\partial t},\\
&\bar{\bm{Q}}^{n+1}=\bar{\bm{Q}}^n+\Delta t \mcR({\bm{F}}^n)+\frac16 \Delta t^2
\left(\frac{\partial \mcR({\bm{F}}^n)}{\partial t}+2\frac{\partial \mcR({\bm{F}}^*)}{\partial t}\right),
\end{split}
\end{equation}
where ${\bm{F}}^*={\bm{F}}(\bm{Q}^*,t)$ is the flux at the intermediate  stage $t^*=t^n+\Delta{t}/2$ in the time interval $[t^n, t^{n+1}]$. The aforementioned two-stage temporal discretization has been proved to achieve fourth-order time accuracy for the hyperbolic conservation law \cite{JLi201638}. Assuming that the computational mesh does not change with time, $\mcR$ is a linear function of ${\bm{F}}$, giving $\partial \mcR({\bm{F}})/\partial t=\mcR(\partial{\bm{F}}/\partial t)$. Thus, the problem remaining is to determine the flux and its time derivative at each flux point.

Based on the gas-kinetic theory, both the macroscopic conservative variables $\bm{Q}$ and the numerical flux can be obtained from the gas distribution function $f$ with the following relations
\begin{equation}\label{eq_f_Q_F}
\begin{split}
  \bm{Q}=\int f \bm{\psi}\mathrm{d}\Xi,~~~
  {\bm{F}_m}=\int u_m f\bm{\psi}\mathrm{d}\Xi,
\end{split}
\end{equation}
where $f=f(\bm{x},t,\bm{u},\bm{\xi})$ is a function of physical space $\bm{x}$, time $t$, particle velocity $\bm{u}$, and the internal DoFs $\bm{\xi}$, $\bm{\psi}={\left(1,\bm{u},(\bm{u}^2+\bm{\xi}^2)/2\right)}^T$ is the vector of moments, and $\rmd\Xi=\rmd\bm{u}\rmd\bm{\xi}$ is the element of the phase space. For convenience, the summation convention is adopted in this subsection, such as $\bm{x}=(x, y)=(x_1, x_2), \bm{u}=(u,v)=(u_1,u_2)$. The governing equation of the gas distribution function is the 2D BGK equation \cite{Bhatnagar1954}
\begin{equation}\label{eq_bgk}
\frac{\partial f}{\partial t}+u_m \frac{\partial f}{\partial x_m}= \frac{g-f}{\tau},
\end{equation}
where $\tau=\mu/p$ is the collision time related to the viscosity $\mu$ and pressure $p$. The local equilibrium state $g$ corresponds to the macroscopic variables, and can be expressed as
 \begin{equation}
 g=\rho (2\pi RT)^{-(K+2)/2} \rme^{-[ (\bm{u}-\bm{U})^2+\bm{\xi}^2]/(2RT)},
\end{equation}
where $K$ is the total number of $\bm{\xi}$. $R$ is the gas constant. 
The conservation law Eq.(\ref{eq_ConLaw}) can be recovered by taking moments of Eq.(\ref{eq_bgk}), where the collision term $(g-f)/\tau$ vanishes automatically owing to the conservation of mass, moments, and total energy during collisions, i.e., the compatibility condition. In particular, the Naiver-Stokes equations can be recovered through the first-order Chapman-Enskog expansion \cite{KXu200214,KXu2015}
\begin{equation}\label{eq_C_E}
f_{NS}=g-\tau \left(\frac{\partial g}{\partial t}+u_m\frac{\partial g}{\partial x_m} \right).
\end{equation} 
By eliminating the first-order term in Eq.(\ref{eq_C_E}), the Euler equations can be recovered as well. To compute the numerical flux through Eq.(\ref{eq_f_Q_F}), a time-dependent gas distribution function is constructed at each flux point based on the analytical solution of Eq.(\ref{eq_bgk})
\begin{equation}\label{eq_BGK_exact}
\begin{split}
 f(\bm{x},t,\bm{u},\xi)=&\frac{1}{\tau}\int_0^t g(\bm{x}-\bm{u}(t-t'),t',\bm{u},\bm{\xi})\mathrm{e}^{-(t-t')/\tau}\mathrm{d}t'\\
                &+\mathrm{e}^{-t/\tau}f_0(\bm{x}-\bm{u}t,\bm{u},\bm{\xi}),
\end{split}
\end{equation}
where $f_0$ is the initial distribution function at the beginning of each time step ($t=0$), and $g$ is the local equilibrium state. For convenience, the subcell interface is assumed to be perpendicular to the $x$-axis, and the flux point is assumed to be the origin $\bm{x}=0$ in a local coordinate system. Using a first-order Taylor expansion of $f_0$ and $g$ around the flux point and combined with Eq.(\ref{eq_C_E}), the time-dependent gas distribution function can be obtained as follows.
\begin{equation}\label{eq_f_dis}
\begin{split}
 &f(0,t,\bm{u},\bm{\xi})\\
 =&g_0\left(1-\rme^{-t/\tau}
  +((t+\tau)\rme^{-t/\tau}-\tau)a_m u_m+(t-\tau+\tau\rme^{-t/\tau})A\right) \\
        &+\mathrm{e}^{-t/\tau}g_R\left(1-(\tau+t)a_m^R u_m-\tau A^R\right) \mathrm{H}(u) \\
        &+\mathrm{e}^{-t/\tau}g_L\left(1-(\tau+t)a_m^L u_m-\tau A^L\right) \left(1-\mathrm{H}(u)\right),
\end{split}
\end{equation}
where the coefficients $a_ m,~A$ are related to the Taylor expansion of the corresponding Maxwellian functions, i.e., the derivatives of $g_0$, and the coefficients $a_m^L,~A^L$ and $a_m^R,~A^R$ correspond to $g_L$ and $g_R$, respectively. $\mathrm{H}(u)$ is the Heaviside function. These coefficients are determined by the reconstructed conservative variables and the slopes, as well as the compatibility condition. For details about the evaluation of these coefficients, we refer to Ref.\onlinecite{CZhang2018}. The gas-kinetic flux solver is intrinsically multidimensional by involving both normal and tangential spatial derivatives in the construction of the gas distribution function.
Besides, compared to the third-order gas-kinetic flux solver, the construction of the above gas distribution function Eq.(\ref{eq_f_dis}) is much simpler, and the computational cost of flux can be reduced significantly. 

Now, the flux and its time derivative are determined as below.
With Eq.(\ref{eq_f_dis}), a time-dependent flux function can be constructed according to Eq.(\ref{eq_f_Q_F}), denoted as $\bm{F}(\bm{Q}^n,t)$, which is a non-linear function of $t$. A simple fitting method is adopted to obtain the approximated flux, and its first-order time derivative \cite{LPan2016II}, where $\bm{F}(\bm{Q}^n,t)$ is approximated by a linear function $\tilde{\bm{F}}(\bm{Q}^n,t)=\bm{F}^n+(t-t^n)\partial_t\bm{F}^n$ within the time interval $[t^n,t^n+\Delta t]$. Denote the time integration of $\bm{F}(\bm{Q}^n,t)$ within $[t_n,t_n+\delta]$ as
\begin{equation}
\label{eq_S2O4_F(t)_int}
\widehat{\bm{F}}(\bm{Q}^n,\delta)=\int_{t^n}^{t^n+\delta}\bm{F}(\bm{Q}^n,t)\rmd{t},
\end{equation}
Then, an equation set can be obtained
\begin{equation}
\label{equ:S2O4_linear_F(t)_int_II}
\begin{split}
 &\frac12\Delta{t}\bm{F}^n+\frac18\Delta{t^2}\partial_t\bm{F}^n=\widehat{\bm{F}}(\bm{Q}^n,\frac12\Delta{t}),\\
 &\Delta{t}\bm{F}^n+\frac12\Delta{t^2}\partial_t\bm{F}^n=\widehat{\bm{F}}(\bm{Q}^n,\Delta t),
\end{split}
\end{equation}
where the left-hand side is the time integration of $\tilde{\bm{F}}(\bm{Q}^n,t)$ within $[t^n,t^n+\Delta t/2]$ and $[t^n,t^n+\Delta t]$, while the right-hand side is obtained according to Eq.(\ref{eq_S2O4_F(t)_int}). By solving the equation set Eq.(\ref{equ:S2O4_linear_F(t)_int_II}), we have
\begin{equation}
\label{equ:S2O4_linear_F(t)_int_III}
\begin{split}
 &\bm{F}^n=\frac{1}{\Delta t}(4\widehat{\bm{F}}(\bm{Q}^n,\frac12\Delta{t})-\widehat{\bm{F}}(\bm{Q}^n,\Delta t)),\\
 &\partial_t\bm{F}^n=\frac{4}{\Delta t^2}(\widehat{\bm{F}}(\bm{Q}^n,\Delta t)-2\widehat{\bm{F}}(\bm{Q}^n,\frac12\Delta{t})).
\end{split}
\end{equation}  
Similarly, the approximated flux $\bm{F}^{\ast}$ and its time derivative $\partial_t\bm{F}^{\ast}$ can be obtained by simply replacing the superscript $n$ in Eq.(\ref{equ:S2O4_linear_F(t)_int_III}) with $\ast$. Finally, the two-stage time stepping method for updating subcell averages can be summarized as 
\begin{equation}\label{eq_S2O4_final}
\begin{split}
&\bar{\bm{Q}}_{i,j}^*=\bar{\bm{Q}}_{i,j}^n+\mcR_{i,j}(\widehat{\bm{F}}),\\
&\bar{\bm{Q}}_{i,j}^{n+1}=\bar{\bm{Q}}_{i,j}^n+\mcR_{i,j}(\widetilde{\bm{F}}),
\end{split}
\end{equation}
where the flux at each flux point is respectively computed according to 
\begin{equation}
\label{equ:S2O4_flux_final}
\begin{split}
  \widehat{\bm{F}}=&\frac12\Delta{t}\bm{F}^n+\frac18\Delta{t^2}\partial_t\bm{F}^n=\widehat{\bm{F}}(\bm{Q}^n,\frac12\Delta{t}),\\
  \widetilde{\bm{F}}=&\frac{8}{3}\widehat{\bm{F}}(\bm{Q}^n,\frac12\Delta{t}) - \frac{1}{3}\widehat{\bm{F}}(\bm{Q}^n,\Delta t) -\frac{8}{3}\widehat{\bm{F}}(\bm{Q}^*,\frac12\Delta{t})\\   &+\frac{4}{3}\widehat{\bm{F}}(\bm{Q}^*,\Delta t).
\end{split}
\end{equation} 
Although there are two stages in the current scheme, the flux evaluation can still be more efficient than the previous single-stage SCFV-GKS because the second-order gas-kinetic flux is much cheaper. Besides, the robustnes of SCFV-GKS can also be enhanced. 

\section{Numerical tests}
In this section, several benchmark cases are tested to validate the performance of the current scheme. For inviscid flows, the collision time $\tau$ is computed by
\begin{equation}\label{eq_tau_inv}
  \tau=\epsilon_1 \Delta t+\epsilon_2\frac{|p^L-p^R|}{|p^L+p^R|}\Delta t,
\end{equation}
where $\epsilon_1$ is set as $1.0$ for $\tau$ included in the exponential term $\mathrm{e}^{-t/\tau}$ in Eq.(\ref{eq_f_dis}) to provide the required numerical dissipation; otherwise, $\epsilon_1$ is set as $10^{-10}$ to better approximate the inviscid assumption. The second term in Eq.(\ref{eq_tau_inv}) represents the artificial viscosity, where $p^L$ and $p^R$ are the pressure at the left and right sides of the subcell interface. The coefficient $\epsilon_2$ is set to 10.0. For viscous flows, the first term $\epsilon_1 \Delta t$ should be replaced by the molecular viscosity $\mu/p$, where $\mu$ is the dynamic viscous coefficient, and $p$ is the pressure at the subcell interface. 

The CFL number is set to $0.2$ for all test cases when the time step $\Delta t$ is computed with reference to the mesh size of main cells, denoted as $h$. Note that it is more reasonable to consider the size of subcells (or CVs) when compared to traditional FV methods. Thus, the corresponding CFL number used in the current scheme is actually $0.4$ with reference to the size of subcells, i.e., $h/2$. According to our tests, the upper limit of the CFL number is empirically 0.35 with the size of main cells (or 0.7 with the size of subcells) for a stable time marching, which is much higher than the DG and CPR methods, and which is comparable with traditional high-order FV methods.
The ratio of specific heat is set to $\gamma=1.4$. 
In the accuracy test, 
 the $L_1$ error and $L_2$ error are computed by
\begin{equation}\label{eq_L1L2error}
\begin{split}
  &L_1~\mathrm{error}=\frac{\sum_{i=1}^{N}\sum_{j=1}^{4}|\bar{q}_{i,j}-\bar{q}_{i,j}^e|}{4N},\\
  &L_2~\mathrm{error}=\sqrt{\frac{\sum_{i=1}^{N}\sum_{j=1}^{4}(\bar{q}_{i,j}-\bar{q}_{i,j}^e)^2}{4N}},
\end{split}
\end{equation}
where $\bar{q}_{i,j}$ and $\bar{q}_{i,j}^e$ indicate the averaged numerical solution and the analytical solution on subcell $\Omega_{i,j}$, respectively, and $N$ is the number of main cells. The two-dimensional (2D) contours of flow fields are plotted based on subcell averaged solutions.

\begin{figure}
  \centering
  \includegraphics[scale=0.05]{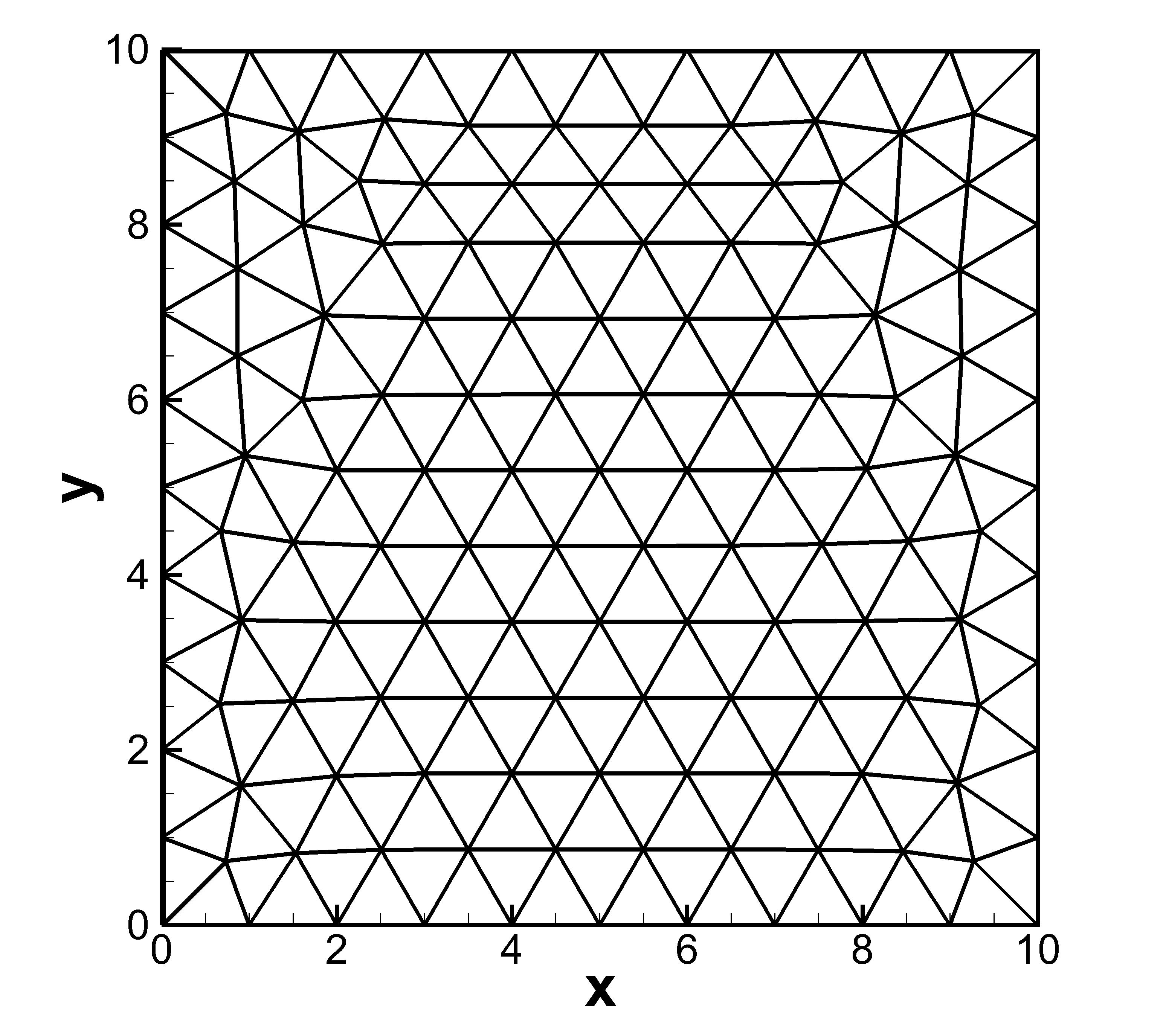}\\
  \caption{Sample mesh for the isentropic vortex propagation.}
  \label{fig:mesh_vortex}
\end{figure}

\subsection{Isentropic vortex propagation}
The isentropic vortex propagation \cite{QWShu1998} was simulated to verify the accuracy and efficiency of the current scheme.
An isentropic vortex was added to the mean flow $(\rho,U,V,p)=(1,1,1,1)$ by introducing perturbations in the velocity $U,~V$ and temperature $T=p/\rho$, but without perturbations in the entropy $S=p/{\rho^{\gamma}}$, i.e.,
\begin{equation}\label{eq_vortex}
\begin{split}
    &(\delta U,\delta V)=\frac{\tilde{\varepsilon}}{2\pi}\mathrm{e}^{(1-r^2)/2}(-\hat{y},\hat{x}),\\
     &\delta T=-\frac{(\gamma-1)\tilde{\varepsilon}^2}{8\gamma\pi^2}\mathrm{e}^{1-r^2},\\
     &\delta S=0,
\end{split}
\end{equation}
where the vortex strength $\tilde{\varepsilon}=5$, $r^2=\hat{x}^2+\hat{y}^2,~\hat{x}=x-5,~\hat{y}=y-5$. The exact solution is that the vortex propagates with constant velocity $(U,V)=(1,1)$. The computational domain is $[0,10]\times[0,10]$, and the periodic boundary condition is adopted for all boundaries. A sample mesh is shown in Figure \ref{fig:mesh_vortex}. The mesh is refined by splitting each cell into four similar finer cells. The computational time is $t=10$. Note that for such a smooth flow problem, no troubled cells are marked by the shock detector, and the HR limiter is not activated. Nevertheless, we also test this case by applying the limiter on all cells artificially to validate its accuracy. 
The density errors and convergence orders are presented in Table~\ref{Accuracy_vortex} and~\ref{Accuracy_vortex_limiter}. 
 The expected order of accuracy is achieved by the current scheme whether or not the limiter is applied. With the limiter, the errors only slightly increase, which means that the high accuracy can still be maintained in smooth regions even if there are cells that are incorrectly marked as troubled cells for shock capturing problems. Thus, the dependence on the shock detector can be reduced effectively.

\begin{table}
\caption{\label{Accuracy_vortex}Accuracy test in the isentropic vortex propagation.}
\renewcommand\arraystretch{1.3}
\begin{ruledtabular}
\begin{tabular}{ccccc}
h & $L_1$ error & Order & $L_2$ error & Order\\
\hline
0.5                & 7.05E-05 &       & 1.20E-04 &      \\
0.25               & 3.95E-06 & 4.16  & 7.39E-06 & 4.02 \\
0.125              & 2.69E-07 & 3.87  & 5.23E-07 & 3.82 \\
0.0625             & 1.87E-08 & 3.84  & 3.64E-08 & 3.85  \\
\end{tabular}
\end{ruledtabular}
\end{table}

\begin{table}
\caption{\label{Accuracy_vortex_limiter}Accuracy test with limiter in the isentropic vortex propagation.}
\renewcommand\arraystretch{1.3}
\begin{ruledtabular}
\begin{tabular}{ccccc}
h & $L_1$ error & Order & $L_2$ error & Order\\
\hline
0.5          & 1.27E-04 &       & 2.80E-04 &       \\
0.25         & 6.84E-06 & 4.22  & 1.88E-05 & 3.89  \\
0.125        & 3.37E-07 & 4.34  & 9.66E-07 & 4.29  \\
0.0625       & 2.47E-08 & 3.77  & 5.07E-08 & 4.25  \\ 
\end{tabular}
\end{ruledtabular}
\end{table}

\begin{figure}
  \centering
  \includegraphics[scale=0.35]{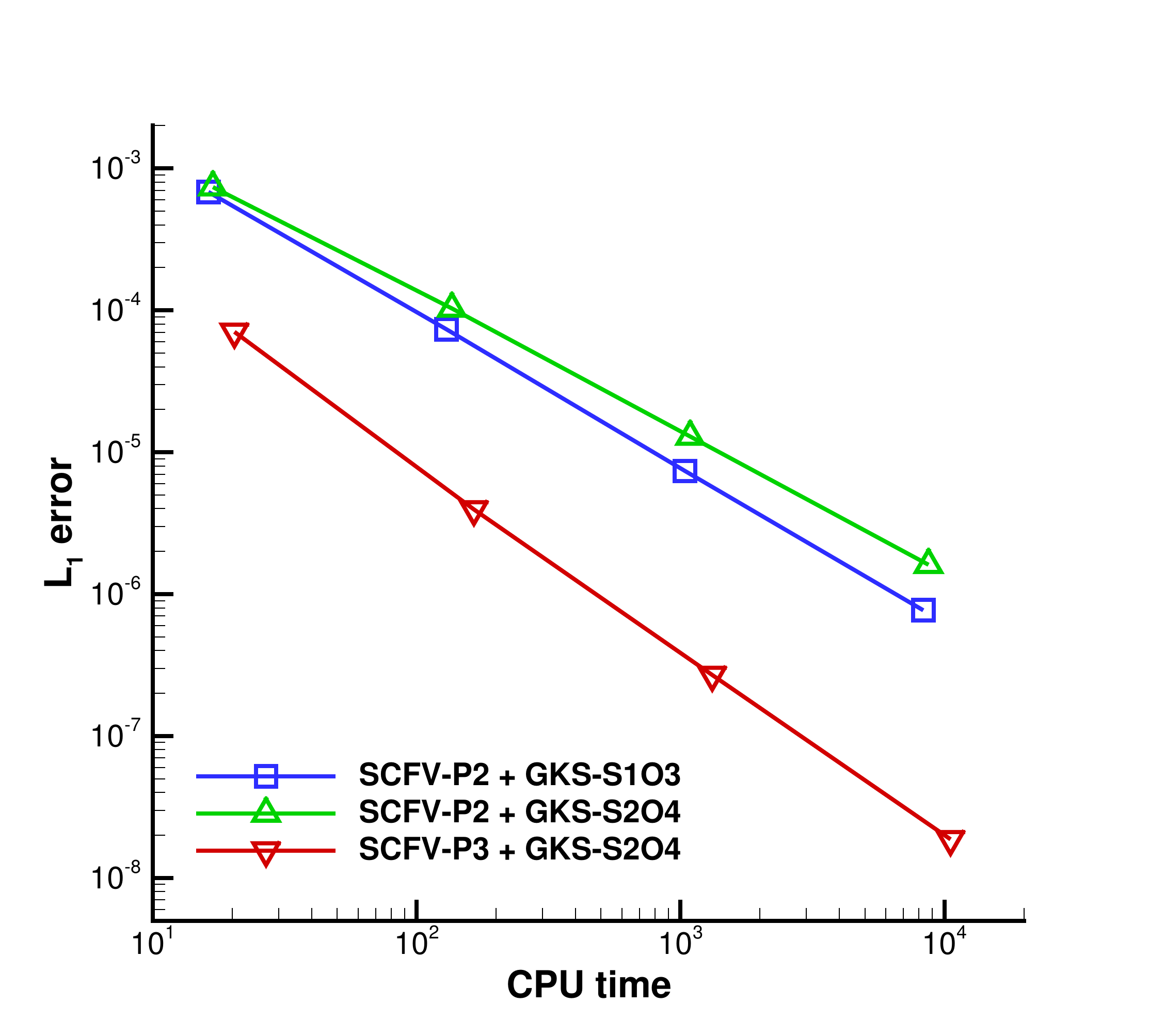}\\
  \caption{Error vs. CPU time in the isentropic vortex propagation.}
  \label{fig:CPU_vortex}
\end{figure}

To estimate the efficiency of the current two-stage SCFV-GKS, Figure \ref{fig:CPU_vortex} shows the relations between the error and the CPU time. For comparison, the result obtained by the single-stage SCFV-GKS is also presented, in which the CLS reconstruction is adopted. Because the computational cost introduced by the limiter has been reduced as much as possible by using the shock detector, for convenience, only the results without the limiter are considered here. It can be seen that the two-stage fourth-order SCFV-GKS is much more efficient than the two kinds of third-order SCFV-GKS. The superiority becomes more distinct as the error decreases. Besides, the single-stage third-order SCFV-GKS is more efficient than the two-stage third-order SCFV-GKS due to the use of the CLS reconstruction, rather than the WLS reconstruction adopted in this study.
With the CLS reconstruction, the computational cost of flux can also be reduced effectively for the single-stage SCFV-GKS. Nevertheless, in the current study, we focus on developing the temporal-spatial fourth-order SCFV-GKS, which is much more efficient. Owing to the two-stage temporal discretization, a simpler and more efficient way of achieving fourth-order accuracy is constructed for the SCFV method. More detailed investigations about the influence of different reconstruction techniques are still needed in the future.

\subsection{Tirarev-Toro problem}
\begin{figure*} 
  \centering
  \begin{varwidth}[t]{\textwidth}
  \vspace{0pt}
  \includegraphics[scale=0.3]{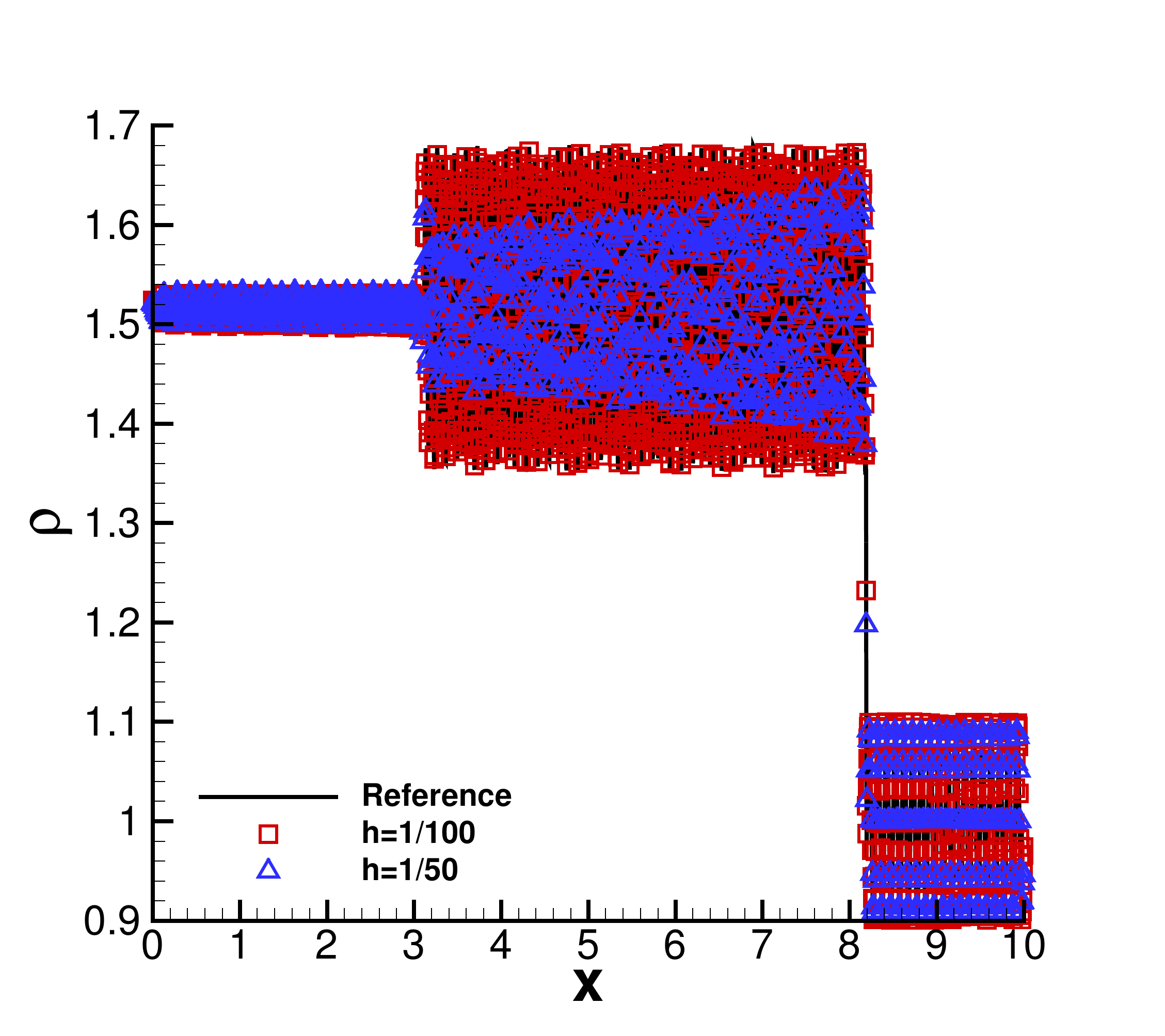}\\
  \end{varwidth}
  \qquad
  \begin{varwidth}[t]{\textwidth}
  \vspace{0pt}
  \includegraphics[scale=0.3]{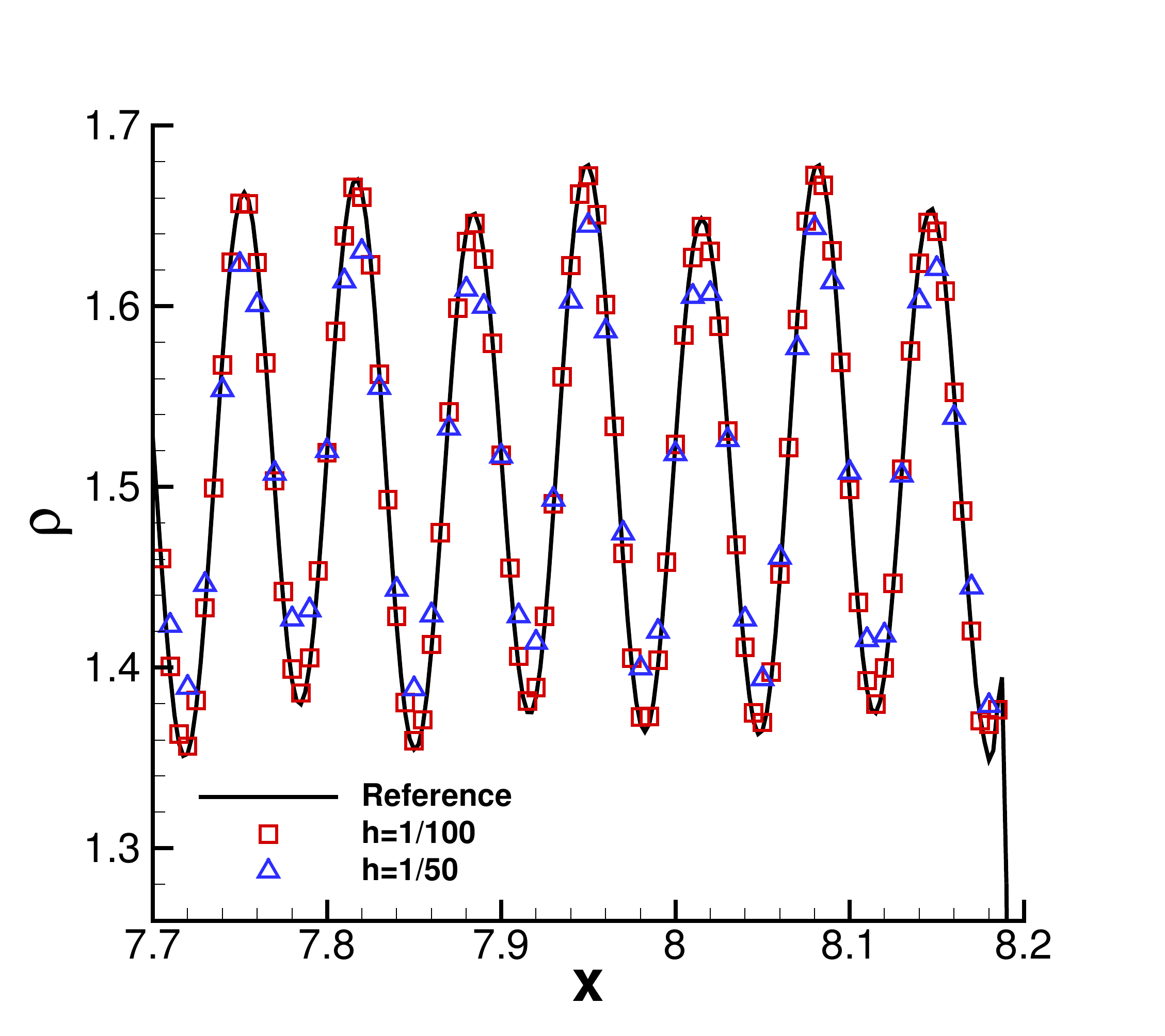}\\
  \end{varwidth}
  \caption{Density distribution along the horizontal centerline at t=5 in the Tirarev-Toro problem.}
  \label{fig:Toro_rho_1D}
\end{figure*}

To confirm the ability of the current scheme to capture high-frequency waves, the Tirarev-Toro problem \cite{VTitarev2014} was tested. The initial condition is
\begin{equation}\label{eq_Tirarev_Toro_initial}
(\rho,U,V,p)=\begin{cases}
(1.515695,0.523346,0,1.805),& 0\leq x \leq 0.5, \\
(1+0.1\sin(20\pi x),0,0,1),& 0.5< x \leq 10.
\end{cases}
\end{equation}
The computational domain is $[0,10]\times[0,0.1]$ with the mesh obtained using the simple triangulation of a rectangular mesh. Figure \ref{fig:Toro_rho_1D} presents the density distribution at $t=5.0$ along the horizontal centerline with the mesh size $h=1/50$ and $1/100$. The reference data is computed by the current scheme with $h=1/1000$. The high-frequency waves are captured by the current scheme accurately, especially with $h=1/100$, which the density profile matches very well with the reference data. With the same mesh size $h=1/100$, the result is slightly better than that obtained by the 6th-order and 8th-order compact GKS \cite{FXZhao2019}. Even the result with $h=1/50$ is comparable with that obtained by the WENO-7JS \cite{FXZhao2019} with $h=1/100$. The results demonstrate the high resolution of the current scheme.

\subsection{Blast wave problem}
\begin{figure*}
  \centering
  \begin{varwidth}[t]{\textwidth}
  \vspace{0pt}
  \includegraphics[scale=0.3]{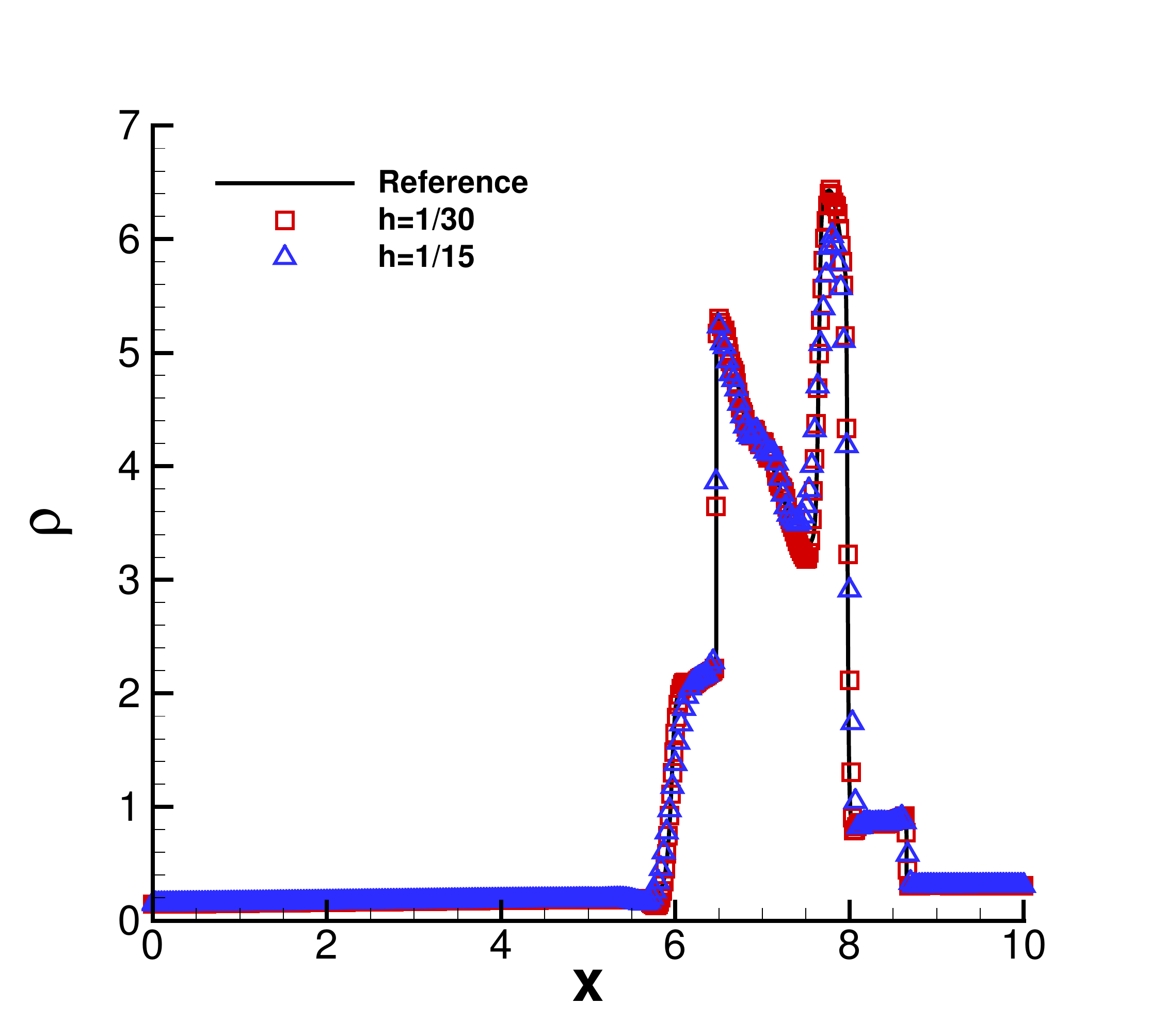}\\
  \end{varwidth}
  \qquad
  \begin{varwidth}[t]{\textwidth}
  \vspace{0pt}
  \includegraphics[scale=0.3]{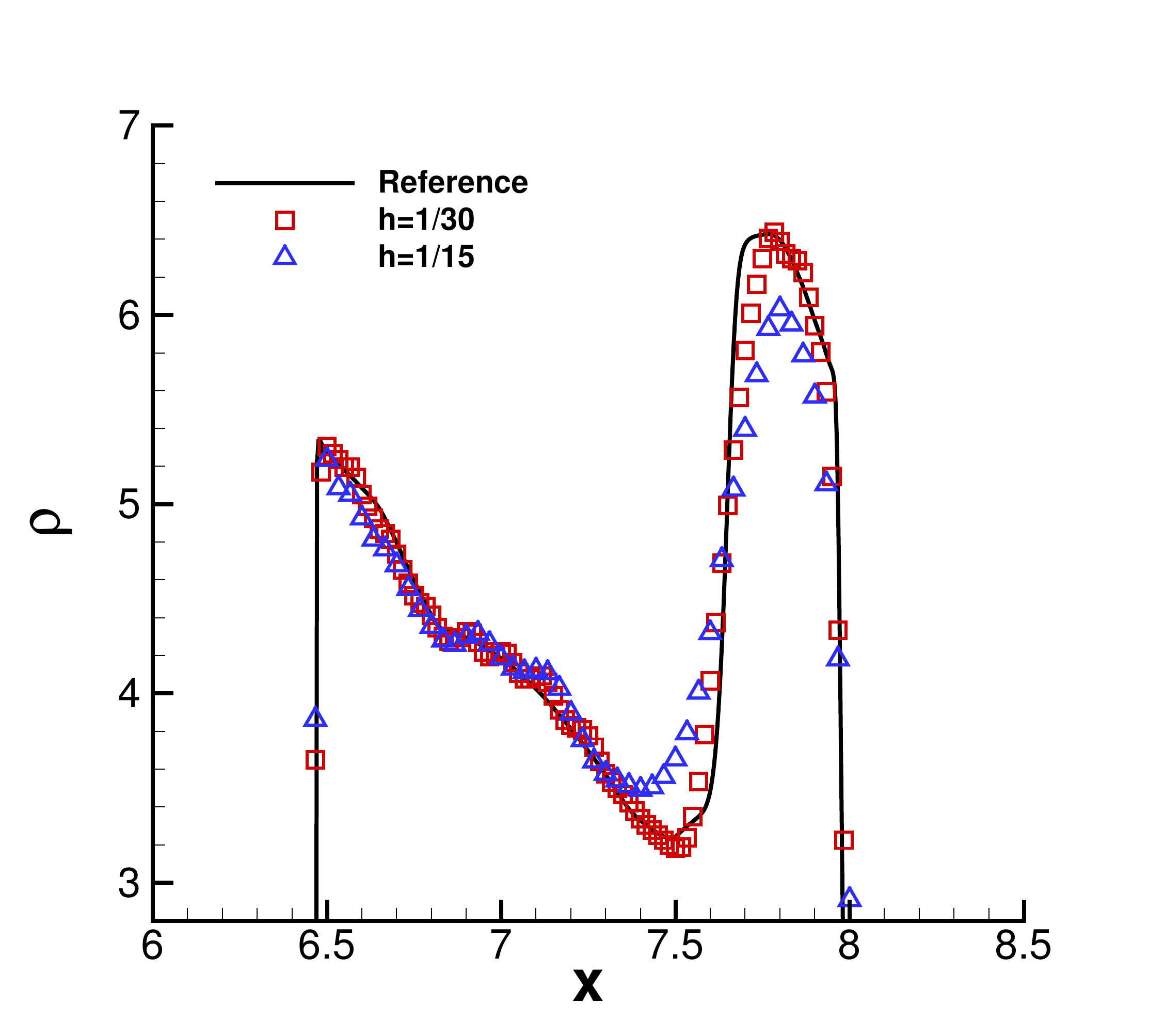}\\
  \end{varwidth}
  \caption{Density distribution along the horizontal centerline at t=0.38 in the blast wave problem.}
  \label{fig:blast_rho_U_1D}
\end{figure*}

Next, the blast wave problem \cite{PWoodward1984} was tested to validate the robustness of numerical schemes for capturing extremely strong shock waves. The initial condition is
\begin{equation}\label{eq_Blast_initial}
(\rho,U,V,p)=\begin{cases}
(1,0,0,1000),& 0\leq x \leq 1, \\
(1,0,0,0.01),& 1< x \leq 9, \\
(1,0,0,100),& 9< x \leq 10.
\end{cases}
\end{equation}
The computational domain is $[0,10]\times[0,1]$ with the mesh obtained by a simple triangulation of a rectangular mesh. The periodic boundary condition is applied to upper and lower boundaries, and the reflecting boundary condition is applied to the left and right boundaries. The reference data are computed by the one-dimensional high-order GKS (HGKS) \cite{QBLi2010} with the mesh size $h=1/10000$. Figure \ref{fig:blast_rho_U_1D} shows the density distribution at $t=0.38$ along the horizontal centerline $y=0.5$ with the mesh size $h=1/15$ and $h=1/30$. The interaction of strong shock waves is captured by the current scheme with no oscillations. The results with $h=1/30$ agree well with the reference data. 

\begin{figure}
  \centering
  \begin{varwidth}[t]{\textwidth}
  \vspace{0pt}
  \includegraphics[scale=0.1]{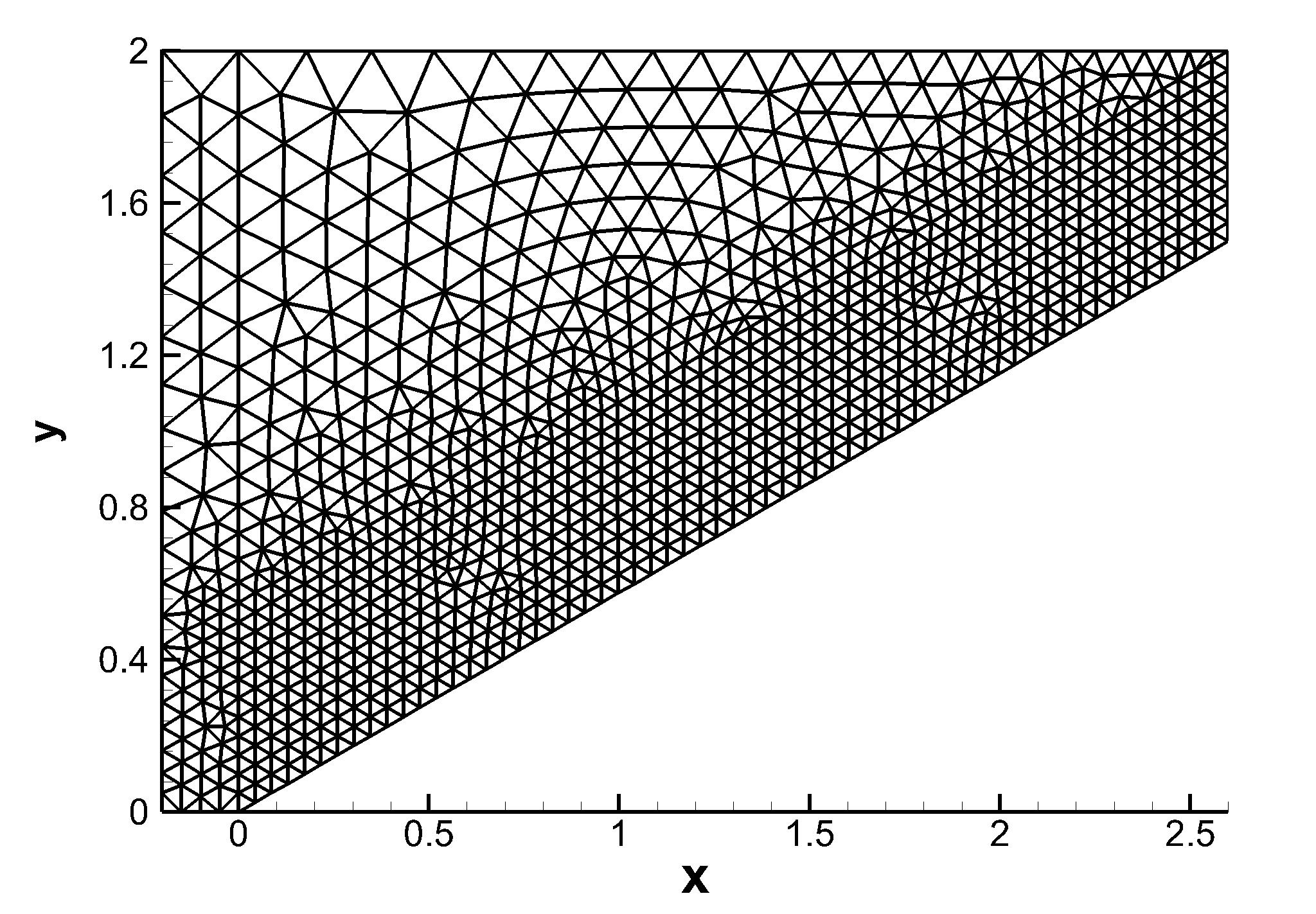}\\
  \end{varwidth}
  \caption{Sample mesh for the double Mach reflection.}
  \label{fig:DMR_sample}
\end{figure}

\subsection{Double Mach reflection}
\begin{figure*}
  \centering
  \begin{varwidth}[t]{\textwidth}
  \vspace{0pt}
  \includegraphics[scale=0.275]{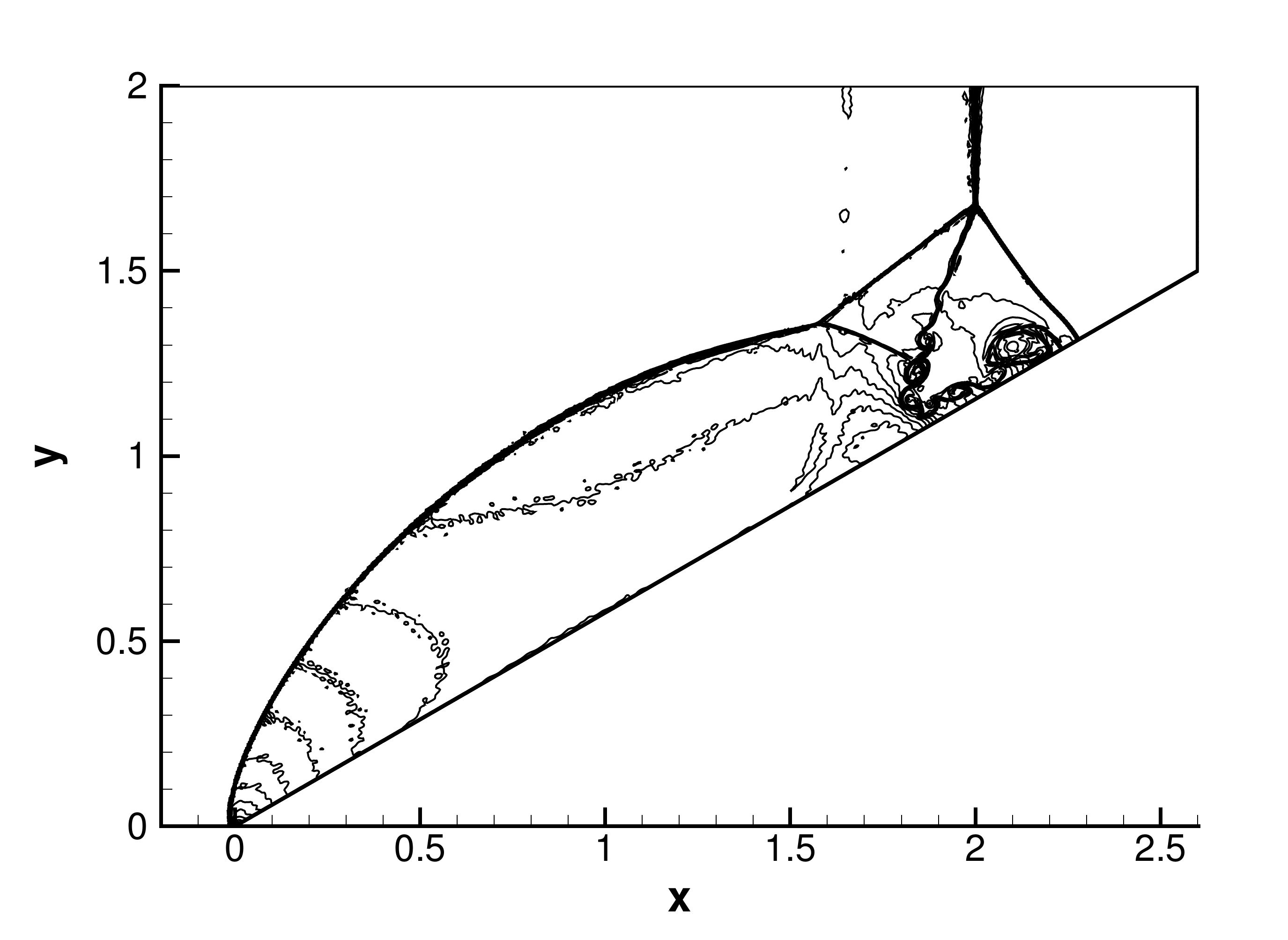}\\
  \end{varwidth}
  \qquad
  \begin{varwidth}[t]{\textwidth}
  \vspace{0pt}
  \includegraphics[scale=0.275]{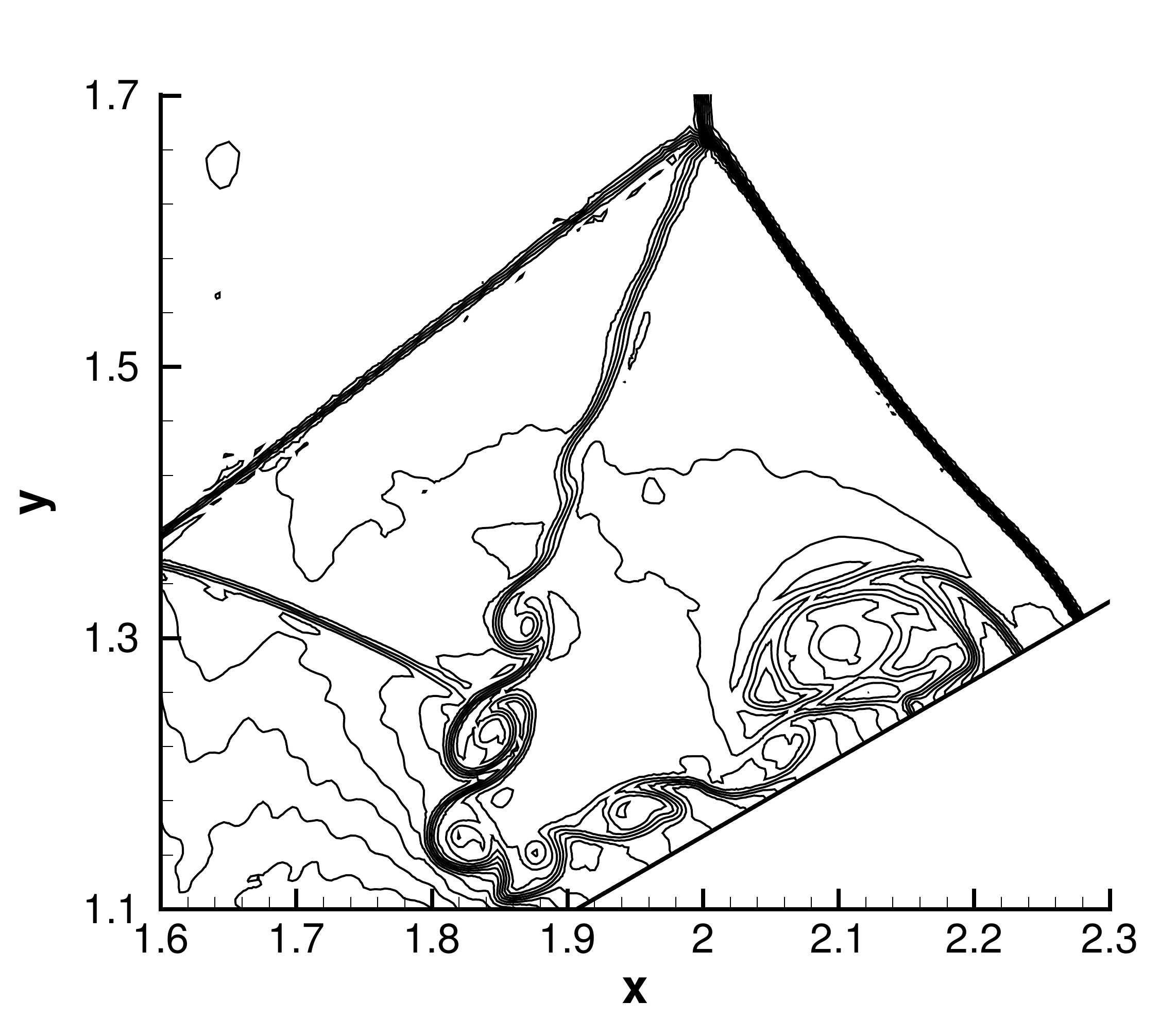}\\
  \end{varwidth}
  \caption{Density contours at $t=0.2$ and the enlarged view near the Mach stem in the double Mach reflection with the mesh size $h=1/160$. Thirty contours are drawn from 2.0 to 22.5.}
  \label{fig:DMR_rho_contours}
\end{figure*}

\begin{figure}
  \centering
  \begin{varwidth}[t]{\textwidth}
  \vspace{0pt}
  \includegraphics[scale=0.23]{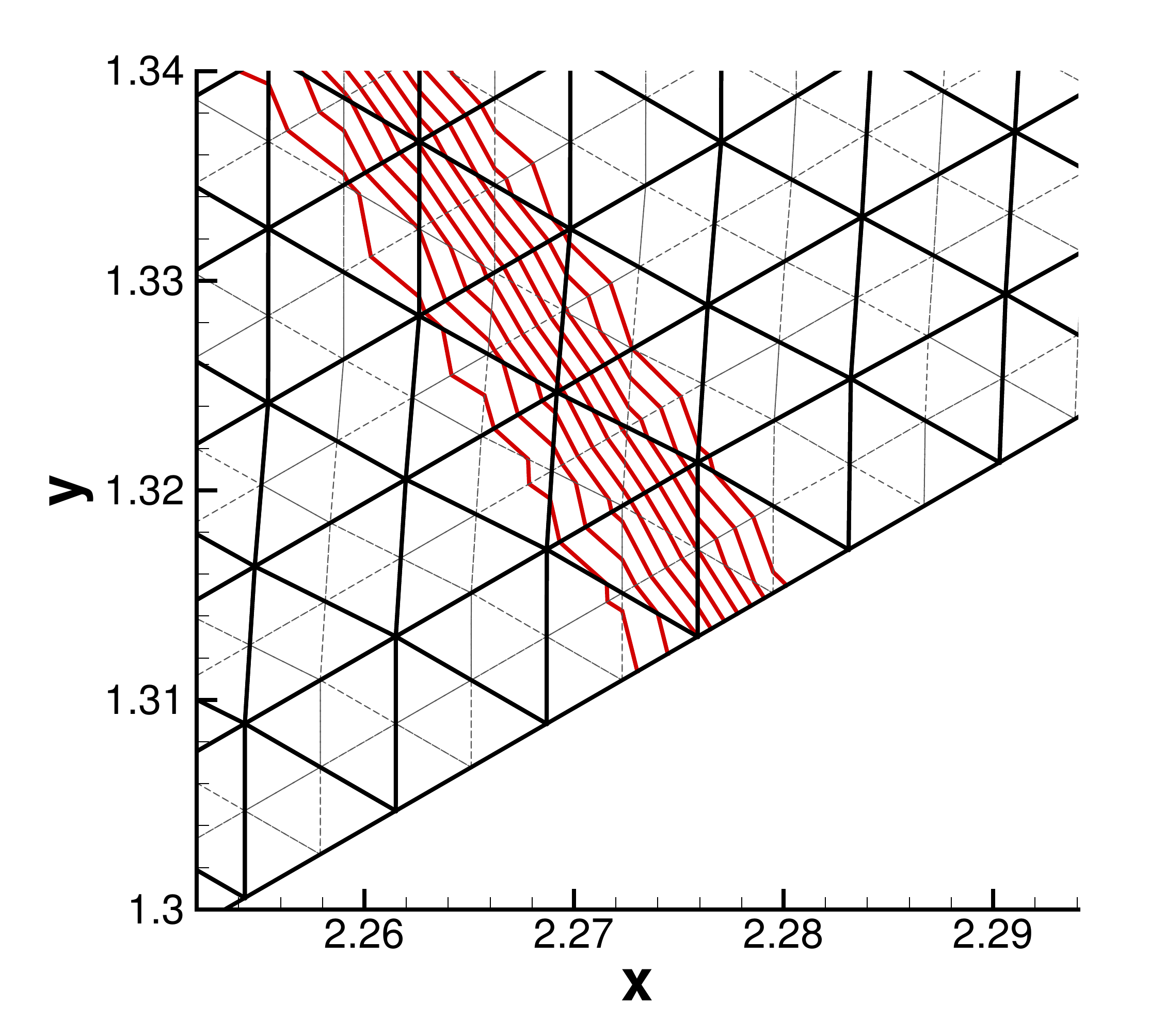}\\
  \end{varwidth}
  \caption{The density contours near the Mach stem in the double Mach reflection.}
  \label{fig:DMR_thickness}
\end{figure}

The double Mach reflection is a 2D benchmark case that is used to validate the robustness and resolution of numerical schemes for shock capturing \cite{PWoodward1984}. A right-moving shock wave with the Mach number $\mathrm{Ma}=10$, initially located at $x=0$,
\begin{equation}\label{eq_DMR_initial}
(\rho,U,V,p)=\begin{cases}
(1.4,0,0,1),& x \leq 0, \\
(8,8.25,0,116.5),& x> 0.
\end{cases}
\end{equation}
which impinges on a $30^{\circ}$ wedge and leads to the double Mach reflection. Figure \ref{fig:DMR_sample} shows the computational domain along with a sample mesh with $h=1/20$, where $h$ indicates the size of main cells uniformly distributed in the region near the wedge, and coarser cells are set elsewhere. The reflecting boundary condition is applied to the wedge and upper boundary. The exact post-shock condition is applied to the left boundary and the bottom boundary from $x=-0.2$ to $x=0$. 
The density contours $t=0.2$ are shown in Figure \ref{fig:DMR_rho_contours} with the mesh size $h=1/120$. The shock waves are captured sharply and the slip line is resolved with high resolution. Besides, Figure  \ref{fig:DMR_thickness} shows the density contours near the Mach stem, where the solid black lines indicate main cells, and the dashed grey lines indicate subcells. It can be observed that the thickness of the shock wave has nearly the same size as the main cell. Because discontinuities can exist inside main cells, the subcell resolution is fully maintained for shock capturing. In contrast, for DG and CPR, the distributions inside main cells need to remain continuous, and the thickness of the shock waves is usually larger than the size of main cells. It is difficult to achieve the high resolution for shock waves as in the SCFV method.

\begin{figure*}
  \centering
  \begin{varwidth}{\textwidth}
  \vspace{0pt}
  \includegraphics[scale=0.23]{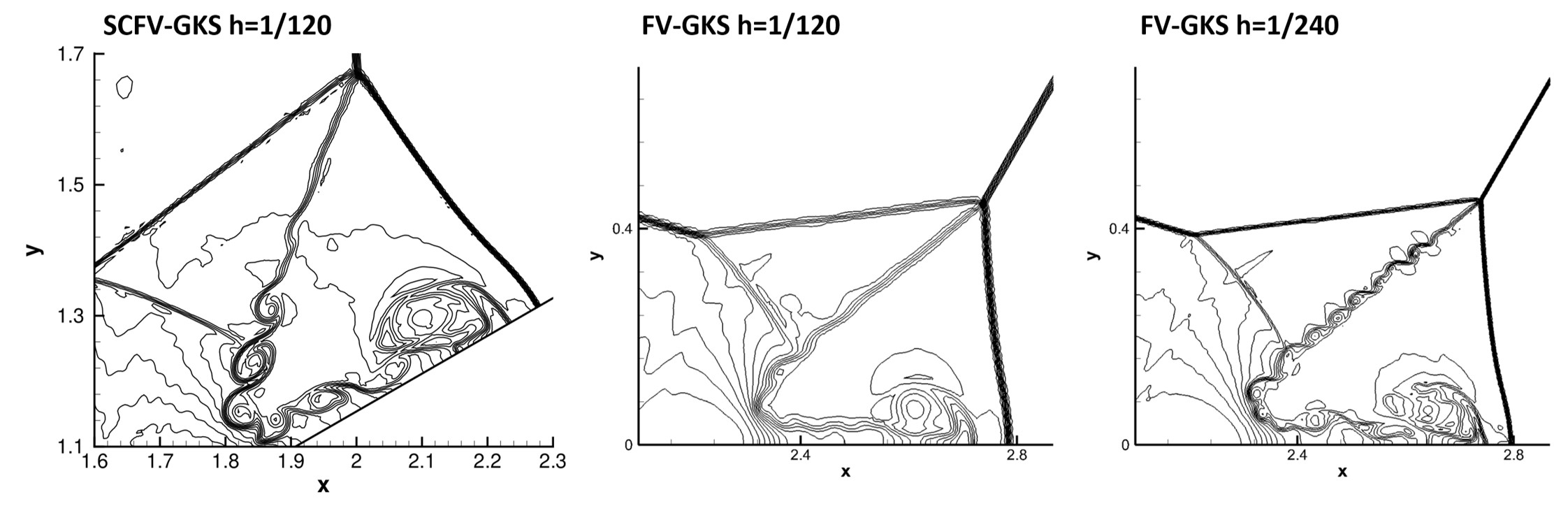}\\
  \end{varwidth}
  \caption{Comparison of the density contours among SCFV-GKS with $h=1/120$ (left), FV-GKS with $h=1/120$ (middle), and $h=1/240$ (right).}
  \label{fig:DMR_rho_contours_compare}
\end{figure*}

In addition, for comparison with traditional FV-GKS, a recently proposed fourth-order FV-GKS based on WENO reconstruction is considered \cite{FXZhao2021Preprint}. As shown in Figure \ref{fig:DMR_rho_contours_compare}, with the same mesh size $h=1/120$, the result obtained by the current scheme is much more accurate than that obtained by FV-GKS owing to the subcell resolution of the current scheme. With a mesh size of $h=1/120$, the size of CVs, i.e., subcells, can be considered to be $1/240$. Therefore, the result is also compared with that obtained by FV-GKS with $h=1/240$. It can be observed that the vortex in the slip line captured by SCFV-GKS has a larger scale than that captured by FV-GKS. Note that only the averaged solutions are updated on each subcell (CV) in the current scheme, while both the averaged solutions and its slopes are updated on each CV in the FV-GKS, which means that the number of DOFs is twice that for FV-GKS, which can help to improve the resolution. However, in FV-GKS, the reconstruction needs to be implemented for each CV, while in the current scheme, a common reconstruction is implemented for four CVs (or subcells), which is more efficient. Besides, only face neighbors are involved in the current scheme for the fourth-order reconstruction, while the stencil in FV-GKS is beyond face neighbors. Therefore, the current scheme is more compact than FV-GKS. 

\subsection{Viscous shock tube flow}
\begin{figure}
  \centering
  \begin{varwidth}[t]{\textwidth}
  \vspace{0pt}
  \includegraphics[scale=0.26]{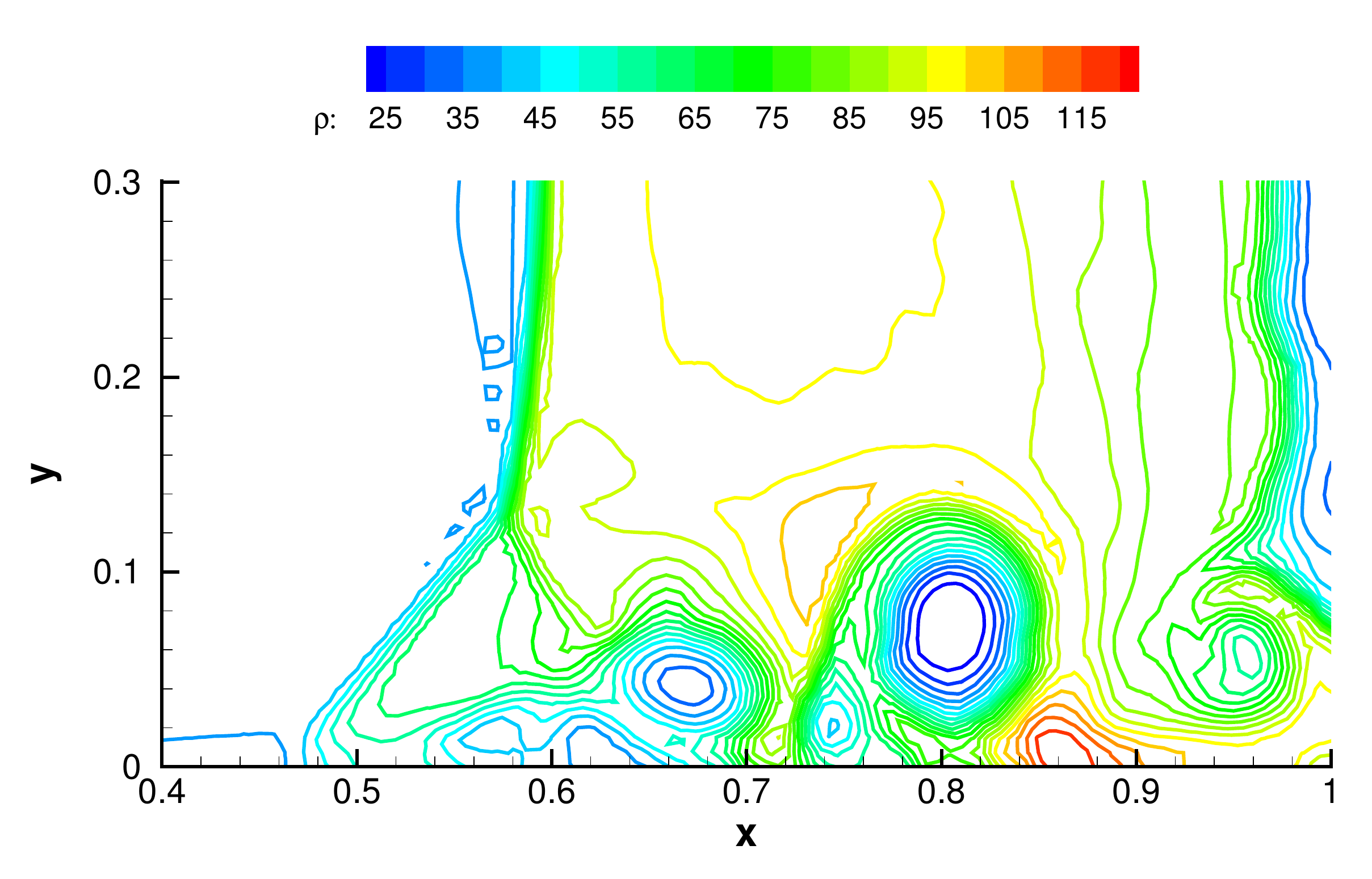}\\
  \end{varwidth}
  \qquad
  \begin{varwidth}[t]{\textwidth}
  \vspace{0pt}
  \includegraphics[scale=0.26]{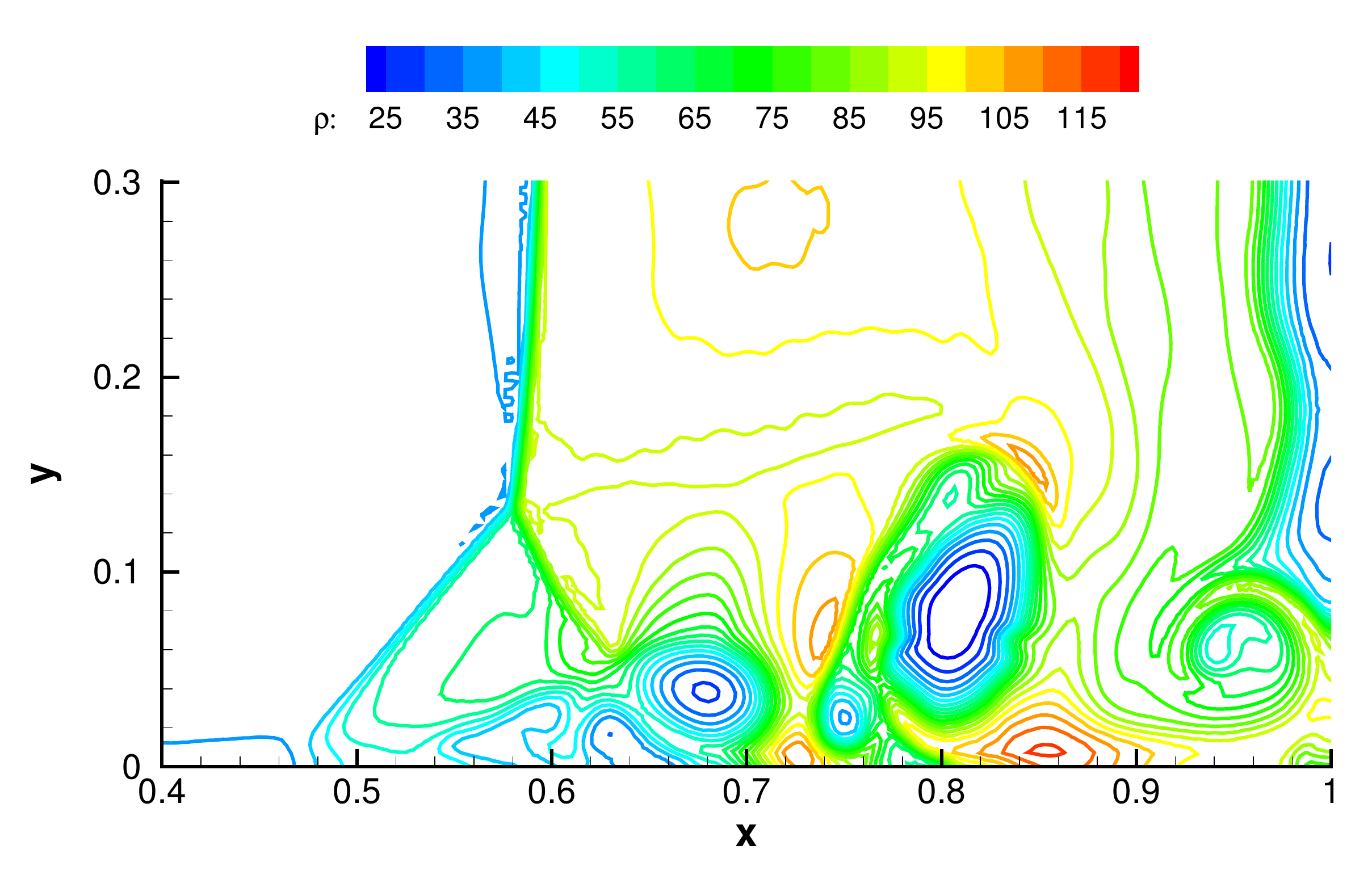}\\
  \end{varwidth}
  \caption{Density contours at t=1 in the viscous shock tube flow with the mesh size $h=1/60$ (top) and $h=1/120$ (bottom), 20 uniform contours from 25 to 120.}
  \label{fig:VST_rho_2D}
\end{figure}

\begin{figure}
  \includegraphics[scale=0.33]{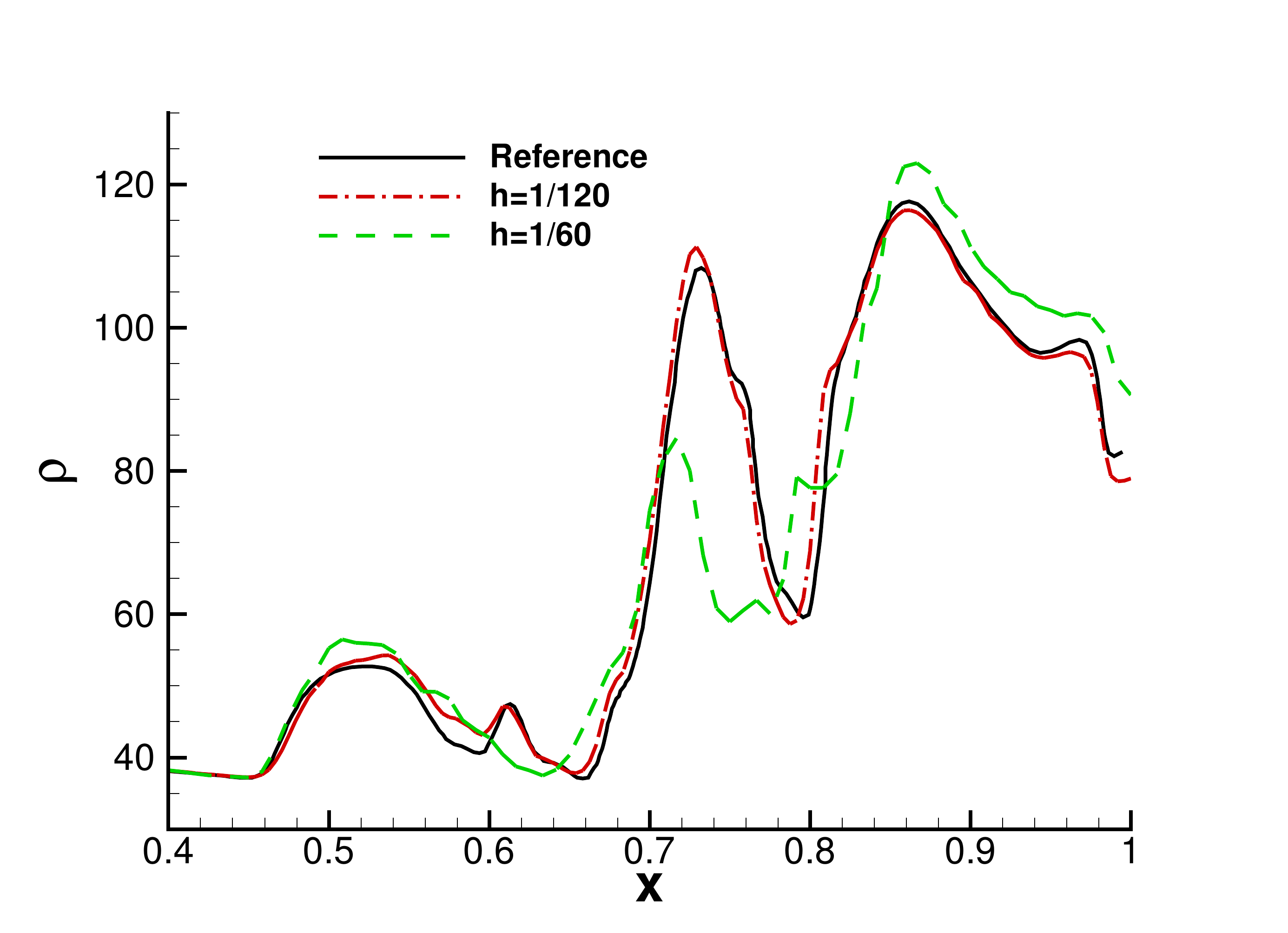}\\
  \caption{Density distribution along the bottom wall in the viscous shock tube flow.}
  \label{fig:VST_rho_wall}
\end{figure}

The viscous shock tube problem is a benchmark case for supersonic viscous flows \cite{Daru2004}, where the flow is bounded by a unit square, and complex unsteady interactions occur between the shock wave and the boundary layer, which requires numerical schemes with both a strong robustness and high resolution. The Reynolds number is set as $\mathrm{Re}=200$ with a constant dynamic viscosity $\mu=0.005$, and the Prandtl number is $\mathrm{Pr}=0.73$. The computational domain is set as $[0,1]\times[0,0.5]$. The symmetrical condition is applied for the upper boundary while the non-slip and adiabatic conditions are applied for other boundaries. 
The initial condition is
\begin{equation}\label{eq_VST_initial}
(\rho,U,V,p)=\begin{cases}
(120,0,0,120/\gamma),& 0\leq x \leq 0.5, \\
(1.2,0,0,1.2/\gamma),& 0.5\leq x \leq 1.
\end{cases}
\end{equation}
The density contours at $t=1$ are presented in Figure \ref{fig:VST_rho_2D} with the mesh size $h=1/60$ and $h=1/120$. The complex flow structures, including the lambda shock and the vortex configurations, are well captured by the current scheme. However, the result with $h=1/60$ is still too dissipative, especially for capturing the primary vortex. By reducing the mesh size to $h=1/120$, the resolution is improved significantly. For additional verification, Figure \ref{fig:VST_rho_wall} shows the density profile along the bottom wall. Table~\ref{VST_Height} presents the estimation of the height of the primary vortex. The results are compared with the reference data \cite{GZZhou2018} provided by a HGKS with $h=1/1500$. As can be seen, the result with $h=1/60$ clearly deviates from the reference data in terms of both the density profile and the height of the primary vortex. In contrast, the result with $h=1/120$ is much more accurate and agrees well with the reference data. Moreover, this flow case is also computed by the compact fourth-order CLSFV and VFV methods \cite{QWang2017} with $h=1/250$. The current result with $h=1/120$ is comparable with that computed by the CLSFV and VFV methods with $h=1/250$. The predicted height of the primary vortex is even closer to the reference data. Considering the advantage of the subcell resolution, the current scheme can achieve accurate solutions with a much coarser mesh compared to traditional FV methods.
The results here demonstrate the good performance of the current scheme in compressible viscous flows. 

\begin{table}
\caption{\label{VST_Height}Height of the primary vortex in the viscous shock tube flow. }
\renewcommand\arraystretch{1.3}
\begin{ruledtabular}
\begin{tabular}{ccc}
Scheme & h & Height\\
\hline
SCFV-GKS & 1/60  & 0.142\\
SCFV-GKS & 1/120 & 0.165\\
Reference & 1/1500 & 0.166\\
\end{tabular}
\end{ruledtabular}
\end{table}

\section{Conclusions}
A compact two-stage fourth-order gas-kinetic SCFV method is developed to solve compressible flows on triangular meshes. The difficulty of compactness faced by traditional FV methods is overcome with the SCFV method by subdividing each cell into a set of subcells in order to increase the number of DoFs for the representation of solution polynomials. The reconstruction is also more efficient than traditional FV because a set of subcells (CVs) share a common reconstruction.
Compared to the single-stage third-order SCFV-GKS, both the accuracy and efficiency are improved significantly by combining the fourth-order compact reconstruction with the second-order accurate flux evolution. With the two-stage fourth-order temporal discretization, we only need to construct the second-order gas distribution function, where the flux evaluation is simplified significantly and the computational cost of flux is reduced to a great extent. Besides, the robustness of SCFV-GKS is enhanced. More importantly, the two-stage temporal discretization provides a more efficient approach to achieve fourth-order time accuracy. Compared to the fourth-stage Runge--Kutta method, one half of the stages can be saved. With the gas-kinetic flux, there is no need to compute the viscous flux for solving viscous flows. Several benchmark cases are tested to verify the performance of the current scheme in compressible flows. The high accuracy and efficiency are validated. 
For shock capturing, the current scheme suppresses oscillations near shock waves effectively; meanwhile, high accuracy can be remained in smooth flow regions. Because discontinuities can exist inside each cell, the subcell resolution is fully maintained, and thus shock waves can be resolved sharply. Hence, the current two-stage fourth-order SCFV-GKS is very promising for the practical simulation of compressible flows. The extension of the proposed method to other types of 2D and 3D unstructured meshes is under consideration in the near future.

\begin{acknowledgments}
Qibing Li was supported by the National Natural Science Foundation of China (Nos. 11672158, 91852109) and the National Key Basic Research and Development Program (No. 2014CB744100). 
Peng Song was supported by the National Natural Science Foundation of China (No. 12031001) and the CAEP foundation (No. CX20200026). 
Jiequan Li was supported by the National Natural Science Foundation of China (Nos. 11771054, 12072042, 91852207), the Sino-German Research Group Project (No. GZ1465), and the National Key Project (GJXM92579). 
The authors would like to thank the technical support of PARATERA and the “Explorer 100” cluster system of Tsinghua National Laboratory for Information Science and Technology.
\end{acknowledgments}

\section*{DATA AVAILABILITY}
The data that support the findings of this study are available
from the corresponding author upon reasonable request.

\section*{references}

\bibliography{aipsamp}

\end{document}